\documentclass[11pt]{article}
\usepackage{tikz}
\usepackage[all]{xy}
\usepackage{amsfonts}
\usepackage{amsfonts}
\usepackage{amsfonts}
\usepackage{mathrsfs}
\usepackage{bm}
\usepackage{cite}
\usepackage[colorlinks,linkcolor=blue,citecolor=blue]{hyperref}
\usepackage{amssymb,amsmath} 
\usepackage{makecell}
\usepackage{tikz}
 \allowdisplaybreaks
\usetikzlibrary{arrows,shapes,chains}
 \allowdisplaybreaks

\catcode`!=11
\let\!int\int \def\int{\displaystyle\!int}
\let\!lim\lim \def\lim{\displaystyle\!lim}
\let\!sum\sum \def\sum{\displaystyle\!sum}
\let\!sup\sup \def\sup{\displaystyle\!sup}
\let\!inf\inf \def\inf{\displaystyle\!inf}
\let\!cap\cap \def\cap{\displaystyle\!cap}
\let\!max\max \def\max{\displaystyle\!max}
\let\!min\min \def\min{\displaystyle\!min}
\let\!frac\frac \def\frac{\displaystyle\!frac}
\catcode`!=12

\let\oldsection\section
\renewcommand\section{\setcounter{equation}{0}\oldsection}

\allowdisplaybreaks
\def\pf{\it{Proof.}\rm\quad}

\def\N{\mathbb{N}}\def\Z{\mathbb{Z}}

\def\sn{\sum\limits_{k=1}^n}

\def\su{\sum\limits_{n=1}^\infty}
\def\sk{\sum\limits_{k=1}^\infty}
\def\sj{\sum\limits_{j=1}^\infty}

\def\z{\zeta}

\def\a{^{(A)}}
\def\B{^{(B)}}
\def\C{^{(C)}}
\def\ab{^{(AB)}}
\def\abc{^{(ABC)}}

\newtheorem{defn}{Definition}[section]
\newtheorem{thm}{Theorem}[section]
\newtheorem{lem}[thm]{Lemma}
\newtheorem{cor}[thm]{Corollary}

\newtheorem{exa}{Example}[section]

\begin{document}
\title {\bf Extensions of Euler Type Sums and Ramanujan Type Sums}
\author{
{Ce Xu\thanks{Corresponding author. Email: cexu2020@ahnu.edu.cn (C. Xu)}}\\[1mm]
\small $*$ School of Mathematics and Statistics, Anhui Normal University,\\ \small Wuhu 241000, People's
Republic of China\\
}

\date{}
\maketitle \noindent{\bf Abstract} We define a new kind of classical digamma function, and establish its some fundamental identities. Then we apply the formulas obtained, and extend tools developed by Flajolet and Salvy to study more general Euler type sums. The main results of Flajolet and Salvy's paper \cite{FS1998} are the immediate corollaries of main results in this paper. Furthermore, we provide some parameterized extensions of Ramanujan-type identities that involve hyperbolic series. Some interesting new consequences and illustrative examples are considered.
\\[2mm]
\noindent{\bf Keywords}: Euler type sum; Polygamma function; (multiple) Zeta value; Contour integration; Residue theorem; Ramanujan-type identity.

\noindent{\bf AMS Subject Classifications (2020):} 65B10; 11M99; 11M06; 11M32.

\section{Introduction and Notations}

Investigation of Euler sums has a long history. The origin of the study of \emph{linear Euler sums (or double zeta star values)}
$$S_{p,q}:=\su \frac{H_n^{(p)}}{n^q}\quad (p\geq 1,\ q\geq 2)$$
which correspond to the $r=1$ case of (\ref{a0}), goes back to the correspondence of Euler with Goldbach in 1742--1743; see Berndt \cite[p. 253]{B1989} for a discussion. Euler elaborated a method to show that the linear sums $S_{p,q}$ can be evaluated in terms of zeta values in the following cases: $p=1$, $p=q$, $p+q$ odd, $p+q$ even but with the pair $(p,q)$ being restricted to the set $\{(2,4),(4,2)\}$. For more details on linear Euler sums, the readers are referred to \cite{BBG1994,BBG1995,FS1998}. Here $H_n^{(p)}$ stands for the $p$-th generalized harmonic number, which is defined by
\begin{align*}
H_n^{(p)}: = \sum\limits_{k=1}^n {\frac{1}{{{k^p}}}} \quad {\rm and}\quad H_0^{(p)}:=0.
\end{align*}
 If $p>1$, the generalized harmonic number $H^{(p)}_n$ converges to the (Riemann) zeta value $\zeta(p)$:
$$\lim\limits_{n\rightarrow \infty}H^{(p)}_n=\zeta(p).$$
When $p=1$, $H_n^{(1)}\equiv H_k$ is the classical harmonic number. A twin sibling of the harmonic number is called alternating harmonic number defined by
\begin{align*}
{ \bar H}_n^{(p)}: = \sum\limits_{k=1}^n {\frac{(-1)^{k-1}}{{{k^p}}}} \quad {\rm and}\quad {\bar H}_0^{(p)}:=0,
\end{align*}
which was introduced in \cite{FS1998}. When taking the limit $n\rightarrow \infty$ in above, we get the so-called the alternating Riemann zeta value
\begin{align*}
\bar \zeta \left( p \right): = \sum\limits_{k = 1}^\infty  {\frac{{{{\left( { - 1} \right)}^{k - 1}}}}
{{{k^p}}}}\quad (p\geq 1).
\end{align*}

Flajolet and Salvy \cite{FS1998} introduced and studied the following a kind of special Dirichlet series that involve harmonic numbers
\begin{align}\label{a0}
{S_{{\bf p},q}} := \sum\limits_{n = 1}^\infty  {\frac{{H_n^{\left( {{p_1}} \right)}H_n^{\left( {{p_2}} \right)} \cdots H_n^{\left( {{p_r}} \right)}}}
{{{n^q}}}},
\end{align}
we call them \emph{generalized (non-alternating) Euler sums}. Moreover, if $r>1$ in (\ref{a0}), they were called \emph{nonlinear Euler sums}. Here ${\bf p}:=(p_1,p_2,\ldots,p_r)\ (r,p_i\in \N, i=1,2,\ldots,r)$ with $p_1\leq p_2\leq \ldots\leq p_r$ and $q\geq 2$. The quantity $w:={p _1} +  \cdots  + {p _r} + q$ is called the weight and the quantity $r$ is called the degree (order). As usual, repeated summands in partitions are indicated by powers, so that for instance
\[{S_{{1^2}{2^3}4,q}} = {S_{112224,q}} = \sum\limits_{n = 1}^\infty  {\frac{{H_n^2[H^{(2)} _n]^3{H^{(4)} _n}}}{{{n^q}}}}. \]
They considered the contour integration involving classical digamma function and used the residue computations to establish more explicit reductions of generalized Euler sums to Euler sums with lower degree. In particular, they proved the famous theorem: a nonlinear Euler sum $S_{p_1p_2\cdots p_r,q}$ reduces to a combination of sums of lower orders whenever the weight $p_1+p_2+\cdots+p_r+q$ and the order $r$ are of the same parity.  The study of nonlinear Euler sums have attracted a lot of research in the area in the last three decades. Some related results may be seen in the works of \cite{M2014,W2017,Xu2017,XW2018,Z2016} and the references therein.

In addition, Flajolet and Salvy \cite{FS1998} also introduced and studied the three \emph{alternating linear Euler sums},
\begin{align}\label{a1}
{S_{{\bar p},q}} :=\su  \frac{{\bar H}^{(p)}_n}{n^q},\quad {S_{p,{\bar q}}} :=\su  \frac{H^{(p)}_n}{n^q}(-1)^{n-1}\quad \text{and}\quad {S_{{\bar p},{\bar q}}} :=\su  \frac{{\bar H}^{(p)}_n}{n^q}(-1)^{n-1}.
\end{align}
Moreover, they shown that the alternating linear Euler sums ${S_{{\bar p},q}}, {S_{p,{\bar q}}} $ and ${S_{{\bar p},{\bar q}}}$ can be evaluated in terms of polynomials of zeta values. For convenience, in (\ref{1.1}), if replace ``$H^{(p_j)}_n$" by ``${\bar H}^{(p_j)}_n$" in the numerator of the summand, we put a ``bar'' on the top of $p_j$. In particular, we put a bar on top of $q$ if there is a sign $(-1)^{n-1}$ appearing in the
denominator on the right. For example,
\begin{align*}
S_{p_1{\bar p}_2p_3{{\bar p}_4},q}:=\su \frac{H^{(p_1)}_n{\bar H}^{(p_2)}_nH^{(p_3)}_n{\bar H}^{(p_4)}_n}{n^q}\quad {\text{and}}\quad S_{{\bar p}_1{\bar p}_2p_3{{\bar p}_4},{\bar q}}:=\su \frac{{\bar H}^{(p_1)}_n{\bar H}^{(p_2)}_nH^{(p_3)}_n{\bar H}^{(p_4)}_n}{n^q}(-1)^{n-1}.
\end{align*}
The sums of types above (one of more the $p_j$ or $q$ barred) are called the \emph{alternating Euler sums}. There are many other researches on alternating Euler sums. For example, Zhao \cite{Z2019} gave the explicit evaluations of all 89 alternating Euler sums with weight $\leq 5$. The author and Wang \cite{XW2018} have developed the Maple package to evaluate the non-alternating Euler sums of weight $2\leq w\leq 16$ and the alternating Euler sums of weight $1\leq w\leq 6$.

In this paper, we define a new kind of classical digamma function. Then, we consider the contour integration involving new digamma function to obtain some new identities of Euler type sums.

Next, we give three definitions. Let $A:=\{a_k\}\ (-\infty < k < \infty)$ be a sequence of complex numbers with ${a_k} = o\left( {{k^\alpha }} \right)\ (\alpha  < 1)$ if $k\rightarrow \pm \infty$. For convenience, we let $A_1$ and $A_2$ to denote the constant sequence $\{(1)^k\}$ and alternating sequence $\{(-1)^k\}$, respectively.

\begin{defn}\label{def1} With $A$ defined above, we define the parametric digamma function $\Psi \left( { - s;A} \right)$ by
\begin{align}\label{1.1}
\Psi \left( { - s;A} \right):= \frac{{{a_0}}}{s} + \sum\limits_{k = 1}^\infty  {\left( {\frac{{{a_k}}}{k} - \frac{{{a_k}}}{{k - s}}} \right)}\quad (s\in\mathbb{ C}\setminus(\N\cup \{0\})).
\end{align}
\end{defn}

If $a_n=1$, then the parametric digamma function $\Psi \left( { - s;A} \right)$ reduces the classical digamma function $\psi \left( { - s} \right)+\gamma$ which is defined as the logarithmic derivative of the well know gamma function
\[\psi \left( s\right) := \frac{d}{{ds}}\left( {\log \Gamma \left( s \right)} \right) = \frac{{\Gamma '\left( s \right)}}{{\Gamma \left( s \right)}}= - \gamma  - \frac{1}{s} + \sum\limits_{k = 1}^\infty  {\left( {\frac{1}{k} - \frac{1}{{k +s}}} \right)} .\]
and it satisfies the complement formula
\[\psi \left( s \right) - \psi \left( { - s} \right) =  - \frac{1}{s} - \pi \cot \left( {\pi s} \right),\]
as well as an expansion at $s=0$ that involves the zeta values:
\[\psi \left( s \right) + \gamma  =  - \frac{1}{s} + \zeta \left( 2 \right)s - \zeta \left( 3 \right){s^2} +  \cdots .\]
From the definition of Riemann zeta function and Hurwitz zeta function, we know that
\[{\psi ^{\left( n \right)}}\left( 1 \right) = {\left( { - 1} \right)^{n + 1}}n!\zeta \left( {n + 1} \right)\quad{\rm and}\quad {\psi ^{\left( n \right)}}\left( z \right) = {\left( { - 1} \right)^{n + 1}}n!\zeta \left( {n + 1;z} \right).\]
The Riemann zeta function and Hurwitz zeta function are defined by
\[\zeta(s):=\sum\limits_{n = 1}^\infty {\frac {1}{n^{s}}}\quad(\Re(s)>1)\]
and
\[\zeta \left( {s;\alpha  + 1} \right): = \sum\limits_{n = 1}^\infty  {\frac{1}{{{{\left( {n + \alpha } \right)}^s}}}} \quad \left( {{\mathop{\Re}\nolimits} \left( s \right) > 1,\ \alpha  \notin \N^ - :=\{-1,-2,\ldots\}} \right).\]
The evaluation of the polygamma function $\psi^{(n)}\left(p/q\right)$ at rational values of the argument
can be explicitly done via a formula as given by K$\ddot{\rm o}$lbig \cite{K1996}, or Andrews, Askey and Roy in
terms of the polylogarithmic or other special functions. Some specific values are listed
in the books \cite{A2000}.

\begin{defn}\label{def2} Define the cotangent function with sequence A by
\begin{align}\label{1.2}
 \pi \cot \left( {\pi s;A} \right) &=  - \frac{{{a_0}}}{s} + \Psi \left( { - s;A} \right) - \Psi \left( {s;A} \right)\nonumber \\
  &= \frac{{{a_0}}}{s} - 2s\sum\limits_{k = 1}^\infty  {\frac{{{a_k}}}{{{k^2} - {s^2}}}}.
\end{align}
\end{defn}
It is clear that if letting $A=A_1$ and $A_2$ in (\ref{1.2}), respectively, then
\begin{align*}
&\cot \left( {\pi s;A_1} \right) = \cot \left( {\pi s} \right)\quad\text{and}\quad \cot \left( {\pi s;A_2} \right) = \csc \left( {\pi s} \right).
\end{align*}

We now provide notations that will be used throughout this paper.

\begin{defn} For nonnegative integers $j\geq 1$ and $n$, we define
\begin{align*}
&D\a(j):=\sk \frac{a_k}{k^j},\quad D\a(1):=0,\\
&E\a_n(j):=\sn \frac{a_{n-k}}{k^j},\quad E\a_0(j):=0,\\
&{\bar E}\a_n(j):=\sn \frac{a_{k-n-1}}{k^j},\quad {\bar E}\a_0(j):=0,\\
&F\a_n(j)= \left\{ {\begin{array}{*{20}{c}} \sk \frac{a_{k+n}-a_k}{k}
   {,\ j=1,}  \\
   {\sk \frac{a_{k+n}}{k^j},\ \ \ \ \ \;\;\;j>1,}  \\
\end{array} } \right.\\
&{\bar F}\a_n(j)= \left\{ {\begin{array}{*{20}{c}} \sk \frac{a_{k-n}-a_k}{k}
   {,\ j=1,}  \\
   {\sk \frac{a_{k-n}}{k^j},\ \ \ \ \ \;\;\;j>1,}  \\
\end{array} } \right.\\
&G\a_n(j):=E\a_n(j)-{\bar E}\a_{n-1}(j)-\frac{a_0}{n^j},\quad G\a_0(j):=0,\\
&L\a_n(j):=F\a_n(j)+(-1)^j{\bar F}\a_n(j),\\
&M\a_n(j):=E\a_n(j)+(-1)^j F\a_n(j),\\
&{\bar M}\a_n(j):={\bar F}\a_n(j)-{\bar E}\a_{n-1}(j),\quad n\geq 1,\\
&R\a_n(j):=G\a_n(j)+(-1)^j L\a_n(j).
\end{align*}
\end{defn}

It is clear that $M\a_n(j)+{\bar M}\a_n(j)=R\a_n(j)+\frac{a_0}{n^j}$, and
if $A=A_1$ and $A_2$, then
\begin{align*}
&M^{(A_1)}_n(j)=H^{(j)}_n+(-1)^j\z(j),\ {\bar M}^{(A_1)}_n(j)=\z(j)-H^{(j)}_{n-1},\ R^{(A_1)}_n(j)=(1+(-1)^j)\z(j),\\
&M^{(A_2)}_n(j)=(-1)^{n-1}{\bar H}^{(j)}_n+(-1)^j \left\{ {\begin{array}{*{20}{c}} (1-(-1)^n)\log(2)
   {,\ j=1,}  \\
   {(-1)^{n-1}{\bar \z}(j),\ \ \ \ \ \;\;\;j>1,}  \\
\end{array} } \right.\\
&{\bar M}^{(A_2)}_n(j)=(-1)^{n}{\bar H}^{(j)}_{n-1}+\left\{ {\begin{array}{*{20}{c}} (1-(-1)^n)\log(2)
   {,\ j=1,}  \\
   {(-1)^{n-1}{\bar \z}(j),\ \ \ \ \ \;\;\;j>1,}  \\
\end{array} } \right.\\
&R^{(A_2)}_n(j)=(-1)^{n-1}(1+(-1)^j){\bar \z}(j).
\end{align*}

In order to distinguish, we let $A^{(l)}:=\{a_k^{(l)}\}\ (-\infty < k < \infty)$ ($l$ is any positive integer) be any sequences of complex numbers with ${a_k^{(l)}} = o\left( {{k^\alpha }} \right)\ (\alpha  < 1)$ if $k\rightarrow \pm \infty$. In above three definitions, the sequence $A$ can be replaced by sequence $A^{(l)}$.

The purposes of this paper are to establish some explicit relations of Euler type sums involving finite sums $M^{(A^{(l_1)})}_n(j),{\bar M}^{(A^{(l_2)})}_n(j)$ and $R^{(A^{(l_3)})}_n(j)$ by using the contour integration that involve the parametric digamma function $\Psi \left( { - s;A^{(l)}} \right)$ and residue computations. Then applying the formulas obtained and letting $A^{(l)}=A_1$ or $A_2$, we can obtain many evaluations for Euler sums and hyperbolic series.

The remainder of this paper is organized as follows.

In the second section we establish several formulas of $\Psi \left( { - s;A} \right)$ and $\cot \left( {\pi s;A} \right)$. In the third and fourth sections we use contour integration and residue computations with the help of formulas of $\Psi \left( { - s;A} \right)$ and $\cot \left( {\pi s;A} \right)$ established in the second section to obtain some closed form representations of linear and quadratic Euler type sums. In the fifth section we evaluate infinite series involving $\coth \left( {\pi s;A} \right)$ by using residue computations. Moreover, we provide some Ramanujan-type identities that involve hyperbolic series.

\section{Several Identities Among Parametric Digamma Function}

In \cite{FS1998}, Flajolet and Salvy used residue computations on large circular contour and specific functions to obtain more independent relations for Euler sums. These functions are of the form $\xi(s)r(s)$, where $r(s):=1/{s^q}$ and $\xi(s)$ is a product of cotangent (or cosecant) and polygamma function. Hence, they gave the following equivalent formulas of cotangent, cosecant and polygamma function at the poles of $\xi(s)r(s)$,
\begin{align*}
&\pi \cot \left( {\pi s} \right)\mathop  = \limits^{s \to n} \frac{1}{{s - n}} - 2\sum\limits_{k = 1}^\infty  {\zeta \left( {2k} \right){{\left( {s - n} \right)}^{2k - 1}}} , \\
&\frac{\pi }
{{\sin \left( {\pi s} \right)}}\mathop  = \limits^{s \to n} {\left( { - 1} \right)^n}\left( {\frac{1}
{{s - n}} + 2\sum\limits_{k = 1}^\infty  {\bar \zeta \left( {2k} \right){{\left( {s - n} \right)}^{2k - 1}}} } \right),\\
&\psi \left( { - s} \right) + \gamma \mathop  = \limits^{s \to n} \frac{1}{{s - n}} + {H_n} + \sum\limits_{k = 1}^\infty  {\left( {{{\left( { - 1} \right)}^k}{\zeta _n}\left( {k + 1} \right) - \zeta \left( {k + 1} \right)} \right){{\left( {s - n} \right)}^k}} ,\;\;n \ge 0 \\
&\psi \left( { - s} \right) + \gamma \mathop  = \limits^{s \to  - n} {H_{n - 1}} + \sum\limits_{k = 1}^\infty  {\left( {{\zeta _{n - 1}}\left( {k + 1} \right) - \zeta \left( {k + 1} \right)} \right){{\left( {s + n} \right)}^k}} ,\;n > 0\\
&\frac{{{\psi ^{\left( {p - 1} \right)}}\left( { - s} \right)}}{{\left( {p - 1} \right)!}}\mathop  = \limits^{s \to n} \frac{1}{{{{\left( {s - n} \right)}^p}}}\left( {1 + {{\left( { - 1} \right)}^p}\sum\limits_{i \ge p} {\left( {\begin{array}{*{20}{c}}
   {i - 1}  \\
   {p - 1}  \\
\end{array}} \right)\left( {\zeta \left( i \right) + {{\left( { - 1} \right)}^i}{\zeta _n}\left( i \right)} \right){{\left( {s - n} \right)}^i}} } \right),\;n \ge 0,\;p > 1 \\
&\frac{{{\psi ^{\left( {p - 1} \right)}}\left( { - s} \right)}}{{\left( {p - 1} \right)!}}\mathop  = \limits^{s \to  - n} {\left( { - 1} \right)^p}\sum\limits_{i \ge 0} {\left( {\begin{array}{*{20}{c}}
   {p - 1 + i}  \\
   {p - 1}  \\
\end{array}} \right)\left( {\zeta \left( {p + i} \right) - {\zeta _{n - 1}}\left( {p + i} \right)} \right){{\left( {s + n} \right)}^i}} ,\;n > 0,\;p > 1.
\end{align*}
In below, we also consider the $\xi(s)r(s)$ (only replace polygamma $\psi^{(p-1)}(-s)$ by parametric polygamma $\Psi^{(p-1)}(-s;A)$) to establish some independent relations for Euler type sums. Thus, we need to obtain the Laurent expansions for parametric polygamma $\Psi^{(p-1)}(-s;A)$) about $s=n$ ($n$ is a any integer).

In this section, we will establish the explicit formulas of parametric polygamma function $\Psi^{(p-1)} \left( { - s;A} \right)$ in terms of infinite series that involve sums $M\a_n(j)$ and ${\bar M}\a(j)$.
The results in this section are basic tools that will be used throughout this paper.
\begin{thm}\label{thm2.1} Let $p\geq 1$ and $n$ be nonnegative integers, if $|s-n|<1$ with $s\neq n$, then
\begin{align}\label{2.1}
\frac{{{\Psi ^{\left( {p - 1} \right)}}\left( { - s;A} \right)}}{{\left( {p - 1} \right)!}}=\frac{1}{{{{\left( {s - n} \right)}^p}}}\left\{a_n-\sj (-1)^j\binom{j+p-2}{p-1} M\a_n(j+p-1)(s-n)^{j+p-1} \right\}.
\end{align}
\end{thm}
\pf From the definition of $\Psi \left( { - s;A} \right)$, if $|s-n|<1$ with $s\neq n$, then it can be rewritten in the form
\begin{align*}
 \Psi \left( { - s;A} \right)=& \frac{{{a_0}}}{s} + \frac{{{a_n}}}{n} - \frac{{{a_n}}}{{n - s}} + \sum\limits_{k = 1}^{n - 1} {\left( {\frac{{{a_k}}}{k} - \frac{{{a_k}}}{{k - s}}} \right)}  + \sum\limits_{k = n + 1}^\infty  {\left( {\frac{{{a_k}}}{k} - \frac{{{a_k}}}{{k - s}}} \right)} \nonumber \\
  =& \frac{{{a_n}}}{{s - n}} + \sum\limits_{k = 1}^n {\frac{{{a_k}}}{k}}  - \sum\limits_{k = 1}^n {\left( {\frac{{{a_k}}}{k} - \frac{{{a_{n - k}}}}{{k + s - n}}} \right)}  + \sum\limits_{k = 1}^\infty  {\left( {\frac{{{a_k}}}{k} - \frac{{{a_{k + n}}}}{{k + n - s}}} \right)}.
\end{align*}
Using the elementary identity
\[(1-x)^{-1} = \sum\limits_{k = 1}^\infty  {{x^{k - 1}}} \quad (|x|<1),\]
we find that
\begin{align}\label{2.2}
\Psi \left( { - s;A} \right)&=\frac{a_n}{s-n}-\sj \left\{(-1)^j E\a_n(j)+F\a_n(j) \right\}(s-n)^{j-1}\nonumber\\
&=\frac{a_n}{s-n}-\sj (-1)^j M\a_n(j)(s-n)^{j-1}.
\end{align}
Then, differentiating (\ref{2.2}) $p-1$ times with respect to $s$ with an elementary calculation, we may easily deduce the desired result.\hfill$\square$

When $n=0$ in (\ref{2.1}), then for $|s|<1$ with $s\neq 0$, we have
\begin{align}\label{2.3}
\frac{{{\Psi ^{\left( {p - 1} \right)}}\left( { - s;A} \right)}}{{\left( {p - 1} \right)!}}=\frac{{{a_0}}}{{{s^p}}} + {\left( { - 1} \right)^p}\sum\limits_{j = 1}^\infty  {\binom{j+p-2}{j-1} D\a(j+p-1){s^{j - 1}}}.
\end{align}

\begin{thm}\label{thm2.2} Let $p$ and $n$ be positive integers, if $|s+n|<1$, then
\begin{align}\label{2.4}
\frac{{{\Psi ^{\left( {p - 1} \right)}}\left( { - s;A} \right)}}{{\left( {p - 1} \right)!}}=(-1)^p \sj \binom{j+p-2}{p-1} {\bar M}\a_n(j+p-1)(s+n)^{j-1}.
\end{align}
\end{thm}
\pf The proof is similar to the previous proof.
For $|s+n|<1$, by a straightforward calculation we see that
\begin{align*}
 \Psi \left( { - s;A} \right) =& \sum\limits_{k = 1}^{n - 1} {\frac{{{a_k}}}{k}}  + \sum\limits_{k = 1}^{n - 1} {\left( {\frac{{{a_{k - n}}}}{{k - s - n}} - \frac{{{a_k}}}{k}} \right)}  + \sum\limits_{k = 1}^\infty  {\left( {\frac{{{a_k}}}{k} - \frac{{{a_{k - n}}}}{{k - s - n}}} \right)}  \\
  =& \sum\limits_{k = 1}^{n - 1} {\frac{{{a_k}}}{k}}  + \sum\limits_{k = 1}^{n - 1} {\left( {\frac{{{a_{k - n}}}}{k}\frac{1}{{1 - \frac{{s + n}}{k}}} - \frac{{{a_k}}}{k}} \right)}  + \sum\limits_{k = 1}^\infty  {\left( {\frac{{{a_k}}}{k} - \frac{{{a_{k - n}}}}{k}\frac{1}{{1 - \frac{{s + n}}{k}}}} \right)} .
\end{align*}
Formally expand the summands on the right side into geometric series to deduce that
\begin{align}\label{2.5}
\Psi \left( { - s;A} \right)&=\sj \left\{{\bar E}\a_{n-1}(j)-{\bar F}\a_n(j) \right\}(s+n)^{j-1}\nonumber\\
&=-\sj {\bar M}\a_n(j)(s-n)^{j-1}.
\end{align}
If differentiating (\ref{2.5}) $p-1$ times with respect to $s$, we obtain (\ref{2.4}) to complete the proof.\hfill$\square$

\begin{thm}\label{thm2.3} With $\cot(\pi s;A)$ defined above, if $|s-n|<1$ with $s\neq n\ (n\in \Z)$, then
\begin{align}\label{2.6}
\pi \cot(\pi s;A)=\frac{a_{|n|}}{s-n}-\sj (-\sigma_n)^j R\a_{|n|}(j)(s-n)^{j-1},
\end{align}
where $\sigma_n$ is defined by the symbol of $n$, namely,
\begin{align*}
\sigma_n:= \left\{ {\begin{array}{*{20}{c}}\ 1
   {,\  n\geq 0}  \\
   {-1,\ n<0.}  \\
\end{array} } \right.
\end{align*}
\end{thm}
\pf Return to (\ref{1.2}), use (\ref{2.1}) and (\ref{2.4}) to arrive at, if $|s-n|<1$ with $s\neq n$ then
\begin{align*}
&\pi \cot(\pi s;A)=\frac{a_n}{s-n}-\sum_{j=1}^\infty \left\{(-1)^jG\a_n(j+1)+L\a_n(j+1) \right\}(s-n)^{j-1},
\end{align*}
if $|s+n|<1$ with $s\neq -n$ then
\begin{align*}
&\pi \cot(\pi s;A)=\frac{a_n}{s+n}-\sum_{j=1}^\infty \left\{G\a_n(j+1)+(-1)^jL\a_n(j+1) \right\}(s+n)^{j-1}.
\end{align*}
Thus, summing these two contributions yields the desired evaluation. \hfill$\square$
\begin{cor}(\cite{FS1998})\label{cor2.4} Let $n$ be a nonnegative integer, then the following formulas hold:
\begin{align*}
&  \bar \psi \left( { - s} \right)= {\left( { - 1} \right)^n}\left\{ {\frac{1}
{{n - s}} + \sum\limits_{k = 0}^\infty  {\left( {{{\left( { - 1} \right)}^k}\bar H_n^{\left( {k + 1} \right)} - \bar \zeta \left( {k + 1} \right)} \right){{\left( {s - n} \right)}^k}} } \right\} \quad (|s-n|<1,\ s\neq n), \\
 & \bar \psi \left( { - s} \right)= {\left( { - 1} \right)^n}\sum\limits_{k = 0}^\infty  {\left( {\bar H_{n - 1}^{\left( {k + 1} \right)} - \bar \zeta \left( {k + 1} \right)} \right){{\left( {s + n} \right)}^k}}\quad (|s+n|<1,\ n>0),
\end{align*}
where $\bar \psi (s)$ denotes the modified digamma function, which is defined by
\begin{align*}
&  \bar \psi \left( s \right): = \sum\limits_{k = 0}^\infty  {\frac{{{{\left( { - 1} \right)}^k}}}
{{s + k}}}=\frac{1}
{2}\psi \left( {\frac{{s + 1}}
{2}} \right) - \frac{1}
{2}\psi \left( {\frac{s}
{2}} \right) \quad (s\in \mathbb{C}).
\end{align*}
\end{cor}
\pf Corollary \ref{cor2.4} follows immediately from (\ref{2.2}) and (\ref{2.5}) with $A=A_2$. \hfill$\square$

\section{Linear and Quadratic Euler Type Sums}

We first state a lemma that will subsequently be used in our proofs of main results.

Flajolet and Salvy \cite{FS1998} defined a kernel function $\xi \left( s \right)$ by the two requirements: 1. $\xi \left( s \right)$ is meromorphic in the whole complex plane. 2. $\xi \left( s \right)$ satisfies $\xi \left( s \right)=o(s)$ over an infinite collection of circles $\left| s \right| = {\rho _k}$ with ${\rho _k} \to \infty $. Applying these two conditions of kernel
function $\xi \left( s \right)$, Flajolet and Salvy discovered the following residue lemma.
\begin{lem}(\cite{FS1998})\label{lem3.1}
Let $\xi \left( s \right)$ be a kernel function and let $r(s)$ be a rational function which is $O(s^{-2})$ at infinity. Then
\begin{align}\label{3.1}
\sum\limits_{\alpha  \in O} {{\mathop{\rm Res}}{{\left[ {r\left( s \right)\xi \left( s \right)},s = \alpha  \right]}}}  + \sum\limits_{\beta  \in S}  {{\mathop{\rm Res}}{{\left[ {r\left( s \right)\xi \left( s \right)},s = \beta  \right]}}}  = 0,
\end{align}
where $S$ is the set of poles of $r(s)$ and $O$ is the set of poles of $\xi \left( s \right)$ that are not poles $r(s)$ . Here ${\mathop{\rm Re}\nolimits} s{\left[ {r\left( s \right)},s = \alpha \right]} $ denotes the residue of $r(s)$ at $s= \alpha$.
\end{lem}

Flajolet and Salvy proved every linear sum $S_{p,q}$ whose weight $p+q$ is odd is expressible as a polynomial in zeta values by applying the kernel function
\[\frac{1}{2}\pi \cot \left( {\pi s} \right)\frac{{{\psi ^{\left( {p - 1} \right)}}\left( { - s} \right)}}{{\left( {p - 1} \right)!}}\]
to the base function $r(s)=s^{-q}$. Elaborating on Euler's work, Nielsen \cite{N1906} also proved this result by a method based on partial fraction expansions. Let $B:=\{b_k\}\ (-\infty < k < \infty)$ be a sequence of complex numbers with ${b_k} = o\left( {{k^\beta }} \right)\ (\beta  < 1)$ if $k\rightarrow \pm \infty$.
Replacing $\cot(\pi s)\psi^{(p-1)}(-s)$ by $\cot(\pi s;A)\Psi^{(p-1)}(-s;B)$, we can get the following theorem.

\begin{thm}\label{thm3.2} For positive integers $p$ and $q>1$,
\begin{align}\label{3.2}
&(-1)^{p+q}\su \frac{{\bar M}\B_n(p)}{n^q}a_n+\su \frac{M\B_n(p)}{n^q}a_n\nonumber\\
&=(-1)^p \sum_{j=1}^p \binom{p+q-j-1}{q-1} \su \frac{R\a_n(j)}{n^{p+q-j}}b_n\nonumber\\ &\quad +(-1)^p 2\sum_{j=1}^{[q/2]} \binom{p+q-2j-1}{p-1} D\a(2j)D\B(p+q-2j)\nonumber\\
&\quad+b_0(1+(-1)^{p+q})D\a(p+q)-(-1)^pa_0\binom{p+q-1}{q}D\B(p+q)\nonumber\\
&\quad-(-1)^p\binom{p+q-1}{p}D\ab(p+q),
\end{align}
where $D\ab(q)$ is defined by
\[D\ab(q):=\su \frac{a_nb_n}{n^q}.\]
\end{thm}
\pf In the context of this paper, the theorem results from applying the kernel function
\[\pi \cot \left( {\pi s};A \right)\frac{{{\Psi ^{\left( {p - 1} \right)}}\left( { - s;B} \right)}}{{\left( {p - 1} \right)!}}\]
to the base function $r(s)=s^{-q}$. Namely, we need to compute the residue of the function
\[{f_1}( {s;A,B}): = \pi \cot \left( {\pi s} ;A\right)\frac{{{\Psi ^{\left( {p - 1} \right)}}\left( { - s;B} \right)}}{{\left( {p - 1} \right)!{s^q}}}.\]
The only singularities are poles at the integers. At a negative integer $-n$ the pole is simple and the residue is
\[{\rm{Res}}\left[ {{f_1}( {s;A,B} ),s =  - n} \right] =(-1)^{p+q} \frac{{\bar M}\B_n(p)}{n^q}a_n.\]
At a positive integer $n$, the pole has order $p+1$ and the residue is
\begin{align*}
 {\rm Res}[f_1(s;A,B),s=n]&=(-1)^p \binom{p+q-1}{p} \frac{a_nb_n}{n^{p+q}} + \frac{{ M}\B_n(p)}{n^q}a_n\\&\quad-(-1)^pb_n \sum_{j=1}^p \binom{p+q-j-1}{q-1} \frac{R\a_n(j)}{n^{p+q-j}}.
\end{align*}
Finally the residue of the pole of order $p+q+1$ at $0$ is found to be
\begin{align*}
{\rm Res}[f_1(s;A,B),s=0]&=(-1)^p \binom{p+q-1}{q}a_0 D\B(p+q)- b_0(1+(-1)^{p+q})D\a(p+q)\\&\quad-(-1)^p 2\sum_{j=1}^{[q/2]} \binom{p+q-2j-1}{p-1} D\a(2j)D\B(p+q-2j).
\end{align*}
Summing these three contributions yields the statement of the theorem. \hfill$\square$

Putting $p=q=2$, $A=A_1$ and $b_n= 1/{2^{|n|}}$, a simple example is as follows:
\begin{align*}
\sum\limits_{n=1}^\infty \frac 1{n^22^n}\sum\limits_{k=1}^n \frac {2^k}{k^2}=\frac {51}{16}\zeta(4)-3{\rm Li}_4\left(\frac1{2}\right)-\frac 3{4}\zeta(2)\log^2(2)-\frac 1{8}\log^4(2).
\end{align*}
Here ${\rm Li}_s(x)$ denotes the polylogarithm function which is defined by
\[{\rm Li}_s\left( x \right): = \sum\limits_{n = 1}^\infty  {\frac{{{x^n}}}{{{n^s}}}} \quad \left( {{\mathop{\Re}\nolimits} \left( s \right) >1,\ |x|\leq 1 } \right),\]
where if ${\Re}(s)=1$, then $x\in \mathbb{C}\setminus ((-\infty,-1)\cup [1,+\infty))$ and $|x|\leq 1$.

If setting $A,B\in\{A_1,A_2\}$ in Theorem \ref{thm3.2}, then we deduce these well-known results of (alternating) linear Euler sums.
\begin{cor}(\cite{FS1998}) For an odd weight $m=p+q\ (q\geq 2)$, the four linear Euler sums are reducible to zeta values,
\begin{align}
\begin{aligned}\label{e1}
 {S_{p,q}}: =& \sum\limits_{n = 1}^\infty  {\frac{{H_n^{\left( p \right)}}}{{{n^q}}}}  \\
  =& \frac{1}{2}\zeta \left( m \right) + \frac{{1 - {{\left( { - 1} \right)}^p}}}{2}\zeta \left( p \right)\zeta \left( q \right) \\
  &+ {\left( { - 1} \right)^p}\sum\limits_{k = 0}^{\left[ {p/2} \right]} {\left( {\begin{array}{*{20}{c}}
   {m - 2k - 1}  \\
   {q - 1}  \\
\end{array}} \right)\zeta \left( {2k} \right)\zeta \left( {m - 2k} \right)}\\
  &+ {\left( { - 1} \right)^p}\sum\limits_{k = 0}^{\left[ {q/2} \right]} {\left( {\begin{array}{*{20}{c}}
   {m - 2k - 1}  \\
   {p - 1}  \\
\end{array}} \right)\zeta \left( {2k} \right)\zeta \left( {m - 2k} \right)} ,
\end{aligned}\\
\begin{aligned}\label{e2}
 {S_{\bar p,\bar q}}: =& \sum\limits_{n = 1}^\infty  {\frac{{\bar H_n^{\left( p \right)}}}{{{n^q}}}{{\left( { - 1} \right)}^{n - 1}}} \\
  =& \frac{1}{2}\zeta \left( m \right) + \frac{{1 - {{\left( { - 1} \right)}^p}}}{2}\bar \zeta \left( p \right)\bar \zeta \left( q \right)\\
  &- {\left( { - 1} \right)^p}\sum\limits_{k = 0}^{\left[ {p/2} \right]} {\left( {\begin{array}{*{20}{c}}
   {m - 2k - 1}  \\
   {q - 1}  \\
\end{array}} \right)\zeta \left( {2k} \right)\bar \zeta \left( {m - 2k} \right)} \\
  &- {\left( { - 1} \right)^p}\sum\limits_{k = 0}^{\left[ {q/2} \right]} {\left( {\begin{array}{*{20}{c}}
   {m - 2k - 1}  \\
   {p - 1}  \\
\end{array}} \right)\zeta \left( {2k} \right)\bar \zeta \left( {m - 2k} \right)} ,
\end{aligned}\\
\begin{aligned}\label{e3}
{S_{p,\bar q}}: =& \sum\limits_{n = 1}^\infty  {\frac{{H_n^{\left( p \right)}}}{{{n^q}}}{{\left( { - 1} \right)}^{n - 1}}} \\
  =& \frac{1}{2}\bar \zeta \left( m \right) + \frac{{1 - {{\left( { - 1} \right)}^p}}}{2}\zeta \left( p \right)\bar \zeta \left( q \right) \\
 & - {\left( { - 1} \right)^p}\sum\limits_{k = 0}^{\left[ {p/2} \right]} {\left( {\begin{array}{*{20}{c}}
   {m - 2k - 1}  \\
   {q - 1}  \\
\end{array}} \right)\bar \zeta \left( {2k} \right){\bar \zeta} \left( {m - 2k} \right)} \\
  &+ {\left( { - 1} \right)^p}\sum\limits_{k = 0}^{\left[ {q/2} \right]} {\left( {\begin{array}{*{20}{c}}
   {m - 2k - 1}  \\
   {p - 1}  \\
\end{array}} \right)\bar \zeta \left( {2k} \right)\zeta \left( {m - 2k} \right)} ,
\end{aligned}\\
\begin{aligned}\label{e4}
{S_{\bar p,q}}: =& \sum\limits_{n = 1}^\infty  {\frac{{\bar H_n^{\left( p \right)}}}{{{n^q}}}}\\
  =& \frac{1}{2}\bar \zeta \left( m \right) + \frac{{1 - {{\left( { - 1} \right)}^p}}}{2}\bar \zeta \left( p \right)\zeta \left( q \right)\\
 & + {\left( { - 1} \right)^p}\sum\limits_{k = 0}^{\left[ {p/2} \right]} {\left( {\begin{array}{*{20}{c}}
   {m - 2k - 1}  \\
   {q - 1}  \\
\end{array}} \right)\bar \zeta \left( {2k} \right)\zeta \left( {m - 2k} \right)} \\
  &- {\left( { - 1} \right)^p}\sum\limits_{k = 0}^{\left[ {q/2} \right]} {\left( {\begin{array}{*{20}{c}}
   {m - 2k - 1}  \\
   {p - 1}  \\
\end{array}} \right)\bar \zeta \left( {2k} \right)\bar \zeta \left( {m - 2k} \right)} ,
\end{aligned}
\end{align}
where $\zeta \left(1\right)$ should be interpreted as $0$ wherever it occurs, and $\zeta \left(0\right)=-1/2$, ${\bar \zeta}\left(0\right)=1/2$.
\end{cor}

Obviously, when $A,B\in\{A_1,A_2\}$ in Theorem \ref{thm3.2}, we see that for even weights, four modified forms of the identity hold, but without any (alternating) linear Euler sum occurring. This gives back well-known nonlinear relations between (alternating) zeta values at even arguments. Flajolet and Salvy \cite{FS1998} applied the kernels $\left(\psi^{(j)}(-s)\right)^2$ to $s^{-q}$ to yield some further relations of even weights. Please see their article for further reference. Next, in a same way, we establish a `duality' sum formula of Euler type sums.

\begin{thm}\label{thm3.4} For positive integers $m$ and $p$,
\begin{align}\label{3.7}
&(-1)^m \sum_{i+j=m-1,\atop i,j\geq 0} \binom{p+i-1}{i}\binom{q+j-1}{j} \su \frac{M\B_n(p+i)}{n^{q+j}}a_n\nonumber\\
&+(-1)^p \sum_{i+j=p-1,\atop i,j\geq 0} \binom{m+i-1}{i}\binom{q+j-1}{j} \su \frac{M\a_n(m+i)}{n^{q+j}}b_n\nonumber\\
&=(-1)^{p+m-1}\binom{p+q+m-2}{q-1} D\ab(p+q+m-1)\nonumber\\ &\quad+a_0(-1)^p\binom{p+q+m-2}{p-1} D\B(p+q+m-1)\nonumber\\ &\quad+b_0(-1)^m\binom{p+q+m-2}{m-1} D\a(p+q+m-1)
\nonumber\\ &\quad+(-1)^{m+p} \sum_{j_1+j_2=q+1,\atop j_1,j_2\geq 1} \binom{j_1+m-2}{j_1-1}\binom{j_2+p-2}{j_2-1}D\a(j_1+m-1)D\B(j_2+p-1).
\end{align}
\end{thm}
\pf Consider
\[f_2(s;A,B):=\frac{\Psi^{(m-1)}(-s;A)\Psi^{(p-1)}(-s;B)}{(m-1)!(p-1)!s^q},\]
which has poles of order $p+m$ at $s=n$ ($n\in\N$). With the help of Theorem \ref{thm2.1}, we compute the residues
\begin{align*}
{\rm Res}[f_2(s;A,B),s=n]&=(-1)^{p+m-1}\binom{p+q+m-2}{q-1}\frac{a_nb_n}{n^{p+q+m-1}}\\
&\quad-(-1)^m \sum_{i+j=m-1,\atop i,j\geq 0} \binom{p+i-1}{i}\binom{q+j-1}{j} \frac{M\B_n(p+i)}{n^{q+j}}a_n\\
&\quad-(-1)^p \sum_{i+j=p-1,\atop i,j\geq 0} \binom{m+i-1}{i}\binom{q+j-1}{j} \frac{M\a_n(m+i)}{n^{q+j}}b_n.
\end{align*}
Clearly, $f_2(s;A,B)$ also has a pole of order $p+q+m$ at $s=0$. Using (\ref{2.3}), we find that
\begin{align*}
&{\rm Res}[f_2(s;A,B),s=0]\\&=a_0(-1)^p\binom{p+q+m-2}{p-1} D\B(p+q+m-1)\nonumber\\ &\quad+b_0(-1)^m\binom{p+q+m-2}{m-1} D\a(p+q+m-1)
\nonumber\\ &\quad+(-1)^{m+p} \sum_{j_1+j_2=q+1,\atop j_1,j_2\geq 1} \binom{j_1+m-2}{j_1-1}\binom{j_2+p-2}{j_2-1}D\a(j_1+m-1)D\B(j_2+p-1).
\end{align*}
Summing these two contributions, we thus immediately deduce (\ref{3.7}) to complete the proof.\hfill$\square$

Hence, setting $A,B\in\{A_1,A_2\}$ in Theorem \ref{thm3.4} yields many linear relations between (alternating) linear Euler sums and polynomials in zeta values. Some illustrate examples please see \cite{BBG1995,BG1996,FS1998}.

\begin{cor}\label{cor3.5} For integer $q>1$,
\begin{align}
 &\frac{3}{2}\left( {{S_{{1^2},q}} - {S_{2,q}}} \right) = \left( {q + 1} \right){S_{1,q + 1}} - \sum\limits_{\scriptstyle {j_1} + {j_2} = q - 1, \hfill \atop
  \scriptstyle {j_1},{j_2} \ge 1 \hfill} {{S_{1,{j_1} + 1}}\zeta \left( {{j_2} + 1} \right)} , \\
 &{S_{{1^3},q}} - 3{S_{12,q}} = q{S_{{1^{\rm{2}}},q + 1}} - \sum\limits_{\scriptstyle {j_1} + {j_2} =q - 1, \hfill \atop
  \scriptstyle {j_1},{j_2} \ge 1 \hfill} {{S_{1,{j_1} + 1}}{S_{1,{j_2} + 1}}} .
\end{align}
\end{cor}
\pf Corollary \ref{cor3.5} follows immediately from Theorem \ref{thm3.4} by setting $m=p=1$, $(a_k,b_k)=(H_k,1)$ and $(a_k,b_k)=(H_k,H_k)$.  \hfill$\square$

Further, by applying the kernels
$$\pi \cot(\pi s)\frac{\psi^{m-1}(-s)\psi^{p-1}(-s)}{(m-1)!(p-1)!}$$
to $s^{-q}$, Flajolet and Salvy gave the explicit formulas of quadratic Euler sums via linear Euler sums and zeta values, see \cite[Theorem 4.2]{FS1998}.
Now, we evaluate more general relation for quadratic Euler type sums in the same manner as in the above.  Let $C:=\{c_k\}\ (-\infty < k < \infty)$ be a sequence of complex numbers with ${c_k} = o\left( {{k^\lambda }} \right)\ (\lambda < 1)$ if $k\rightarrow \pm \infty$, and let
\[D\abc(q):=\su \frac{a_nb_nc_n}{n^q}.\]

\begin{thm}\label{thm3.6} Let $m,p$ and $q>1$ be positive integers with $A,B$ and $C$ defined above, we have
\begin{align}
&(-1)^{p+q+m} \su \frac{{\bar M}\B_n(m){\bar M}\C_n(p)}{n^q} a_n+\su \frac{M\B_n(m)M\C_n(p)}{n^q}a_n\nonumber\\
&+(-1)^{p+m}\binom{p+q+m-1}{q-1} D\abc(p+q+m)\nonumber\\
&+(-1)^m\sum_{j=1}^{m+1} \binom{j+p-2}{p-1}\binom{m+q-j}{q-1}\su \frac{M\C_n(j+p-1)}{n^{m+q-j+1}}a_nb_n\nonumber\\
&+(-1)^p\sum_{j=1}^{p+1} \binom{j+m-2}{m-1}\binom{p+q-j}{q-1}\su \frac{M\B_n(j+m-1)}{n^{p+q-j+1}}a_nc_n\nonumber\\
&-(-1)^{p+m} \sum_{j=1}^{p+m} \binom{p+q+m-j-1}{q-1} \su \frac{R\a_n(j)}{n^{p+q+m-j}}b_nc_n\nonumber\\
&-(-1)^m \sum_{j_1+j_2\leq m+1,\atop j_1,j_2\geq 1} \binom{m+q-j_1-j_2}{q-1}\binom{j_2+p-2}{p-1} \su \frac{R\a_n(j_1)M\C_n(j_2+p-1)}{n^{m+q-j_1-j_2+1}}b_n\nonumber\\
&-(-1)^p \sum_{j_1+j_2\leq p+1,\atop j_1,j_2\geq 1} \binom{p+q-j_1-j_2}{q-1}\binom{j_2+m-2}{m-1} \su \frac{R\a_n(j_1)M\B_n(j_2+m-1)}{n^{p+q-j_1-j_2+1}}c_n\nonumber\\
&+{\rm Res}[f_3(s;A,B,C),s=0]=0,
\end{align}
where
\begin{align}\label{3.11}
&{\rm Res}[f_3(s;A,B,C),s=0]\nonumber\\
&=(-1)^p a_0b_0\binom{p+q+m-1}{p-1} D\C(p+q+m)\nonumber\\&\quad+(-1)^ma_0c_0\binom{p+q+m-1}{m-1}D\B(p+q+m)\nonumber\\&\quad-b_0c_0(1+(-1)^{p+q+m})D\a(p+q+m)\nonumber\\
&\quad-(-1)^p2b_0\sum_{2j_1+j_2=m+q+1,\atop j_1,j_2\geq 1} \binom{j_2+p-2}{j_2-1}D\a(2j_1)D\C(j_2+p-1)\nonumber\\
&\quad-(-1)^m2c_0\sum_{2j_1+j_2=p+q+1,\atop j_1,j_2\geq 1} \binom{j_2+m-2}{j_2-1}D\a(2j_1)D\B(j_2+m-1)\nonumber\\
&\quad+(-1)^{m+p}a_0 \sum_{j_1+j_2=q+2,\atop j_1,j_2\geq 1} \binom{j_1+m-2}{j_1-1}\binom{j_2+p-2}{j_2-1}D\B(j_1+m-1)D\C(j_2+p-1)\nonumber\\
&\quad-(-1)^{m+p}2\sum_{2j_1+j_2+j_3=q+2,\atop j_1,j_2,j_3\geq 1} \binom{j_2+m-2}{j_2-1}\binom{j_3+p-2}{j_3-1}\nonumber\\&\quad\quad\quad\quad\quad\quad\quad\quad\quad\quad\quad\quad\times D\a(2j_1)D\B(j_2+m-1)D\C(j_3+p-1).
\end{align}
\end{thm}
\pf  Consider
\[f_3(s;A,B,C):=\frac{\pi\cot(\pi s;A)\Psi^{(m-1)}(-s;B)\Psi^{(p-1)}(-s;C)}{(m-1)!(p-1)!s^q},\]
which has simple poles at $s=-n$ ($n\in \N$), with residues
\begin{align*}
{\rm Res}[f_3(s;A,B,C),s=-n]=(-1)^{p+q+m} \frac{{\bar M}\B_n(m){\bar M}\C_n(p)}{n^q} a_n.
\end{align*}
Clearly, $s=n$ ($n\in \N$) are the poles of order $m+p+1$ of $f_3(s;A,B,C)$, using Theorems \ref{thm2.1} and \ref{thm2.3}, we find that
\begin{align*}
&{\rm Res}[f_3(s;A,B,C),s=n]\\&=\frac{{M}\B_n(m){M}\C_n(p)}{n^q} a_n+(-1)^{p+m}\binom{p+q+m-1}{q-1} \frac{a_nb_nc_n}{n^{p+q+m}}\\
&+(-1)^m\sum_{j=1}^{m+1} \binom{j+p-2}{p-1}\binom{m+q-j}{q-1}\frac{M\C_n(j+p-1)}{n^{m+q-j+1}}a_nb_n\\
&+(-1)^p\sum_{j=1}^{p+1} \binom{j+m-2}{m-1}\binom{p+q-j}{q-1}\frac{M\B_n(j+m-1)}{n^{p+q-j+1}}a_nc_n\\
&-(-1)^{p+m} \sum_{j=1}^{p+m} \binom{p+q+m-j-1}{q-1}\frac{R\a_n(j)}{n^{p+q+m-j}}b_nc_n\\
&-(-1)^m \sum_{j_1+j_2\leq m+1,\atop j_1,j_2\geq 1} \binom{m+q-j_1-j_2}{q-1}\binom{j_2+p-2}{p-1}\frac{R\a_n(j_1)M\C_n(j_2+p-1)}{n^{m+q-j_1-j_2+1}}b_n\\
&-(-1)^p \sum_{j_1+j_2\leq p+1,\atop j_1,j_2\geq 1} \binom{p+q-j_1-j_2}{q-1}\binom{j_2+m-2}{m-1}\frac{R\a_n(j_1)M\B_n(j_2+m-1)}{n^{p+q-j_1-j_2+1}}c_n.
\end{align*}
Moreover, $f_3(s;A,B,C)$ also has a pole of order $p+q+m+1$ at $s=0$, using (\ref{2.3}) and Theorem \ref{thm2.3}, we arrive at (\ref{3.11}). Thus, the desired evaluation holds.\hfill$\square$

Hence, if setting $A=B=C=A_1$ in Theorem \ref{thm3.6}, then it becomes the Theorem 4.2 of Flajolet and Salvy \cite{FS1998}. Further, we obtain the following description.
\begin{cor} If $p+q+m$ is even, and $q>1,m,p$ are positive integers, then the (alternating) quadratic Euler sums
\begin{align*}
&\su \frac{H^{(m)}_nH^{(p)}_n}{n^q},\ \su \frac{{\bar H}^{(m)}_nH^{(p)}_n}{n^q},\ \su \frac{{\bar H}^{(m)}_n{\bar H}^{(p)}_n}{n^q}, \\
&\su \frac{H^{(m)}_nH^{(p)}_n}{n^q}(-1)^{n-1},\ \su \frac{{\bar H}^{(m)}_nH^{(p)}_n}{n^q}(-1)^{n-1},\ \su \frac{{\bar H}^{(m)}_n{\bar H}^{(p)}_n}{n^q}(-1)^{n-1}
\end{align*}
are reducible to (alternating) linear Euler sums.
\end{cor}
\pf Letting $A,B,C\in\{A_1,A_2\}$ in Theorem \ref{thm3.6} yields the desired description. \hfill$\square$

It is obvious that (alternating) linear Euler sums reduce to (alternating) zeta values in the case of an odd weight, while (alternating) quadratic Euler sums reduce to (alternating) linear Euler sums in the case
of an even weight. Flajolet and Salvy also shown a result to the effect that such reductions of order are general, but not explicit formulas, for detail see \cite[Theorem 5.3]{FS1998}. Moreover, by using the above expressions and the AMZVs computation function in the Mathematica
package \emph{MultipleZetaValues} developed by Au \cite{Au2020} recently, we obtain the evaluations of more nonlinear Euler sums. For example, we can get the following cases
\begin{align*}
S_{1\bar2,3}&=-\frac{5}{2} \z(\bar5,1)+\frac{13 \zeta (3)^2}{32}+\frac{143 \pi ^6}{181440},\\
S_{1 \bar2,\bar3}&=\frac{1}{3} \pi ^2 \text{Li}_4\left(\frac{1}{2}\right)-\frac{211 \zeta (3)^2}{64}+\frac{7}{24} \pi ^2 \zeta (3) \log (2)\\&\quad+\frac{103 \pi ^6}{45360}+\frac{1}{72} \pi ^2 \log ^4(2)-\frac{1}{72} \pi ^4 \log ^2(2),\\
S_{\bar1 \bar2,3}&=-\frac{7}{2}\z(\bar5,1)+\frac{153 \zeta (3)^2}{64}-\frac{5}{12} \pi ^2 \zeta (3) \log (2)\\&\quad+\frac{217}{32} \zeta (5) \log (2)-\frac{701 \pi ^6}{181440},\\
S_{\bar1 2,\bar3}&=\frac{5}{2} \z(\bar5,1)-\frac{1}{3} \pi ^2 \text{Li}_4\left(\frac{1}{2}\right)-\frac{61 \zeta (3)^2}{64}+\frac{5}{16} \pi ^2 \zeta (3) \log (2)\\&\quad-\frac{155}{32} \zeta (5) \log (2)+\frac{23 \pi ^6}{5184}-\frac{1}{72} \pi ^2 \log ^4(2)+\frac{1}{72} \pi ^4 \log ^2(2),\\
S_{\bar1\bar2,\bar3}&=-5 \z(\bar5,1)-\frac{1}{6} \pi ^2 \text{Li}_4\left(\frac{1}{2}\right)+\frac{39 \zeta (3)^2}{16}-\frac{9}{16} \pi ^2 \zeta (3) \log (2)\\&\quad+\frac{217}{32} \zeta (5) \log (2)-\frac{5 \pi ^6}{2268}-\frac{1}{144} \pi ^2 \log ^4(2)+\frac{1}{144} \pi ^4 \log ^2(2),\\
S_{\bar1 2,3}&=5 \z(\bar5,1)-\frac{7 \zeta (3)^2}{4}+\frac{29}{48} \pi ^2 \zeta (3) \log (2)-\frac{155}{32} \zeta (5) \log (2)+\frac{781 \pi ^6}{362880},\\
S_{1 2,\bar3}&=\frac3{2} \z(\bar5,1)+\frac{1}{6} \pi ^2 \text{Li}_4\left(\frac{1}{2}\right)-\frac{49 \zeta (3)^2}{64}+\frac{7}{48} \pi ^2 \zeta (3) \log (2)+\frac{29 \pi ^6}{181440}\\&\quad+\frac{1}{144} \pi ^2 \log ^4(2)-\frac{1}{144} \pi ^4 \log ^2(2),\\
S_{\bar1\bar3,\bar2}&=5 \z(\bar5,1)+\frac{1}{3} \pi ^2 \text{Li}_4\left(\frac{1}{2}\right)-\frac{75 \zeta (3)^2}{32}+\frac{3}{4} \pi ^2 \zeta (3) \log (2)-\frac{93}{32} \zeta (5) \log (2)\\&\quad-\frac{71 \pi ^6}{72576}+\frac{1}{72} \pi ^2 \log ^4(2)-\frac{1}{72} \pi ^4 \log ^2(2),
\end{align*}
where $\z(\bar5,1)$ is a alternating double zeta values, see (\ref{AMZV-Defn}).

\section{Formulas for General Euler Type Sums}

As with pervious work, in this section, we prove some explicit relations of Euler type sums for arbitrary degree. We begin with some basic notations. For positive integer sequences ${\bf j}:=(j_1,j_2,\ldots,j_r)$ and ${\bf p}:=(p_1,p_2,\ldots,p_r)$, we let
\begin{align*}
&|{\bf j}|:=j_1+j_2+\cdots+j_r,\quad {\bf p}:=p_1+p_2+\cdots+p_r,\\
&C_r({\bf j};{\bf p}):=\binom{j_1+p_1-2}{p_1-1}\binom{j_2+p_2-2}{p_2-1}\cdots\binom{j_r+p_r-2}{p_r-1},\quad C_0({\bf j};{\bf p}):=1.
\end{align*}
Let $\mathcal{S}_r$ be the symmetric group of all the permutations on $r$ symbols. For a permutation $\sigma\in S_r$, we let
\begin{align*}
&|{\bf j}_\sigma|_l:=j_{\sigma(1)}+j_{\sigma(2)}+\cdots+j_{\sigma(l)},\quad |{\bf p}_\sigma|_l:=p_{\sigma(1)}+p_{\sigma(2)}+\cdots+p_{\sigma(l)},\\
&C_l({\bf j}_\sigma;{\bf p}_\sigma):=\binom{j_{\sigma(1)}+p_{\sigma(1)}-2}{p_{\sigma(1)}-1}\cdots\binom{j_{\sigma(l)}+p_{\sigma(l)}-2}{p_{\sigma(l)}-1},\quad C_0({\bf j}_\sigma;{\bf p}_\sigma):=1,
\end{align*}
where $l=1,2,\ldots,r$ and $|{\bf j}_\sigma|_0=|{\bf p}_\sigma|_0:=0$.

We are now ready to state general evaluations in closed form for general Euler type sums.
\begin{defn} With $A^{(l)}$ defined above, we define
\begin{align}
f(s;{\bf A}):=\frac{\Psi^{(p_1-1)}(-s;A^{(1)})\Psi^{(p_2-1)}(-s;A^{(2)})\cdots \Psi^{(p_r-1)}(-s;A^{(r)})}{(p_1-1)!(p_2-1)!\cdots (p_r-1)!}.
\end{align}
\end{defn}

From Theorems \ref{thm2.1} and \ref{thm2.2}, we find that if $|s-n|<1$ with $s\neq n\ (n\geq 0)$, then
\begin{align}\label{4.2}
f(s;{\bf A})=\sum_{l=0}^r\sum_{\sigma\in \mathcal{S}_r} \sum_{j_{\sigma(1)},\ldots,j_{\sigma{(l)}}=1}^\infty &(-1)^{l+|{\bf j}_\sigma|_l}C_l({\bf j}_\sigma;{\bf p}_\sigma)\prod_{i=1}^l M^{(A^{(\sigma(i))})}_n(j_{\sigma(i)}+p_{\sigma(i)}-1)\prod_{k=l+1}^r a^{(\sigma(k))}_n\nonumber\\
&\times (s-n)^{|{\bf j}_\sigma|_l+|{\bf p}_\sigma|_l-l-|{\bf p}|}.
\end{align}
If $|s+n|<1\ (n\geq 1)$, then
\begin{align}\label{4.3}
f(s;{\bf A})=(-1)^{|{\bf p}|} \sum_{j_1,\ldots,j_r=1}^\infty C_r({\bf j};{\bf p})\prod_{i=1}^r {\bar M}^{(A^{(i)})}_n(j_{i}+p_{i}-1)(s+n)^{|{\bf j}|-r}.
\end{align}

\begin{thm}\label{thm4.1} Let $p_1,p_2,\ldots,p_r,q$ and $r$ be positive integers. For a permutation $\sigma\in S_r$, we have
\begin{align}
&\sum_{l=0}^{r-1} \sum_{\sigma\in \mathcal{S}_r}\sum_{|{\bf j}_\sigma|_l\leq w_\sigma({\bf p},l)} (-1)^{|{\bf p}|-|{\bf p}_\sigma|_l+1}C_l({\bf j}_\sigma;{\bf p}_\sigma)\binom{w_\sigma({\bf p},l)+q-|{\bf j}_\sigma|_l-1}{q-1}\nonumber\\ &\quad\quad\quad\quad\quad\quad\quad\quad\times\su \frac{\prod_{i=1}^l M^{(A^{(\sigma(i))})}_n(j_{\sigma(i)}+p_{\sigma(i)}-1)\prod_{k=l+1}^r a^{(\sigma(k))}_n}{n^{w_\sigma({\bf p},l)+q-|{\bf j}_\sigma|_l}}\nonumber\\
&+\sum_{l=1}^{r} \sum_{\sigma\in \mathcal{S}_r}\sum_{|{\bf j}_\sigma|_l=q+ w_\sigma({\bf p},l)} (-1)^{|{\bf p}_\sigma|_l}C_l({\bf j}_\sigma;{\bf p}_\sigma)\prod_{i=1}^l D^{(A^{(\sigma(i))})}(j_{\sigma(i)}+p_{\sigma(i)}-1)\prod_{k=l+1}^r a^{(\sigma(k))}_0\nonumber\\&=0,
\end{align}
where $w_\sigma({\bf p},l):=|{\bf p}|-|{\bf p}_\sigma|_l+l-1$.
\end{thm}
\pf Applying the kernel $f(s;{\bf A})$ to the base function $r(s)=1/s^q$, and using (\ref{4.2}) and (\ref{4.3}), we achieve the desired expansion after a rather tedious computation.\hfill$\square$

If setting $r=3,p_1=p_2=p_3=1$ and $a^{(1)}_k=a^{(2)}_k=a^{(3)}_k=1,(-1)^k,H_k$ in Theorem \ref{thm4.1}, we get the following examples.
\begin{exa}
For integer $q>1$,
\begin{align}\label{4.5}
 {S_{{1^2},q}} - {S_{2,q}}{\rm{ = }}&p{S_{1,q + 1}} - \frac{{\left( {q - 2} \right)\left( {q + 3} \right)}}{6}\zeta \left( {q + 2} \right) + \zeta \left( 2 \right)\zeta \left( q \right) - \sum\limits_{\scriptstyle {j_1} + {j_2} = q, \hfill \atop
  \scriptstyle {j_1},{j_2} \ge 1 \hfill} {\zeta \left( {{j_1} + 1} \right)\zeta \left( {{j_2} + 1} \right)} \nonumber \\
  &+ \frac{1}{3}\sum\limits_{\scriptstyle {j_1} + {j_2} + {j_3} = q - 1, \hfill \atop
  \scriptstyle {j_1},{j_2},{j_3} \ge 1 \hfill} {\zeta \left( {{j_1} + 1} \right)\zeta \left( {{j_2} + 1} \right)\zeta \left( {{j_{_3}} + 1} \right)} ,
\end{align}
\begin{align}\label{4.6}
 {S_{{{\bar 1}^2},\bar q}} =&  - q{S_{\bar 1,{\overline {q + 1}}}} - {S_{\bar 2,\bar q}} - \frac{{\left( {q - 2} \right)\left( {q + 3} \right)}}{6}\bar \zeta \left( {q + 2} \right) + q\log (2)\left( {\zeta \left( {q + 1} \right) + \bar \zeta \left( {q + 1} \right)} \right) \nonumber\\
  &- 2{\log ^2}(2)\left( {\zeta \left(q \right) + \bar \zeta \left( q \right)} \right) + 2\log (2)\left( {{S_{\bar 1,q}} + {S_{\bar 1,\bar q}}} \right) - \frac{1}{2}\zeta \left( 2 \right)\bar \zeta \left( q \right) \nonumber\\
  &+ \frac{1}{3}\sum\limits_{\scriptstyle {j_1} + {j_2} + {j_3} = q - 1, \hfill \atop
  \scriptstyle {j_1},{j_2},{j_3} \ge 1 \hfill} {\bar \zeta \left( {{j_1} + 1} \right)\bar \zeta \left( {{j_2} + 1} \right)\bar \zeta \left( {{j_{_3}} + 1} \right)} \nonumber \\
  &+ \sum\limits_{\scriptstyle {j_1} + {j_2} = q, \hfill \atop
  \scriptstyle {j_1},{j_2} \ge 1 \hfill} {\bar \zeta \left( {{j_1} + 1} \right)\bar \zeta \left( {{j_2} + 1} \right)} ,
\end{align}
\begin{align}\label{4.7}
 &\frac{{q\left( {q + 1} \right)}}{6}{S_{{1^3},q + 2}} - \frac{q}{2}\left( {{S_{{1^4},q + 1}} - 3{S_{{1^2}2,q+ 1}}} \right) - 2\zeta \left( 3 \right){S_{{1^2},q}} - \zeta \left( 2 \right){S_{{1^3},q}} - \frac{5}{2}{S_{{1^3}2,q}} + 2{S_{{1^2}3,q}} \nonumber\\
 & + \frac{1}{4}\left( {{S_{{1^5},q}} + 9{S_{{{12}^2},q}}} \right) - \frac{1}{3}\sum\limits_{\scriptstyle {j_1} + {j_2} + {j_3} = q- 1, \hfill \atop
  \scriptstyle {j_1},{j_2},{j_3} \ge 1 \hfill} {{S_{1,{j_1} + 1}}{S_{1,{j_2} + 1}}{S_{1,{j_3} + 1}}}  = 0.
\end{align}
\end{exa}

\begin{thm}\label{thm4.2} Let $p_1,p_2,\ldots,p_r,q$ and $r$ be positive integers. For a permutation $\sigma\in S_r$, we have
\begin{align}
&(-1)^{|{\bf p}|+q}\su \frac{\prod_{i=1}^r {\bar M}^{(A^{(i)})}_n(p_{i})}{n^q}a_n+\su \frac{\prod_{i=1}^r { M}^{(A^{(i)})}_n(p_{i})}{n^q}a_n\nonumber\\
&+\sum_{l=0}^{r-1} \sum_{\sigma\in \mathcal{S}_r}\sum_{|{\bf j}_\sigma|_l\leq w'_\sigma({\bf p},l)} (-1)^{|{\bf p}|-|{\bf p}_\sigma|_l}C_l({\bf j}_\sigma;{\bf p}_\sigma)\binom{w'_\sigma({\bf p},l)+q-|{\bf j}_\sigma|_l-1}{q-1}\nonumber\\
 &\quad\quad\quad\quad\quad\quad\quad\quad\times\su \frac{\prod_{i=1}^l M^{(A^{(\sigma(i))})}_n(j_{\sigma(i)}+p_{\sigma(i)}-1)\prod_{k=l+1}^r a^{(\sigma(k))}_n}{n^{w'_\sigma({\bf p},l)+q-|{\bf j}_\sigma|_l}}\nonumber\\
 &-\sum_{l=0}^{r-1} \sum_{\sigma\in \mathcal{S}_r}\sum_{|{\bf j}_\sigma|_l+t\leq w'_\sigma({\bf p},l)} (-1)^{|{\bf p}|-|{\bf p}_\sigma|_l}C_l({\bf j}_\sigma;{\bf p}_\sigma)\binom{w'_\sigma({\bf p},l)+q-|{\bf j}_\sigma|_l-t-1}{q-1}\nonumber\\
  &\quad\quad\quad\quad\quad\quad\quad\quad\times\su \frac{R\a_n(t)\prod_{i=1}^l M^{(A^{(\sigma(i))})}_n(j_{\sigma(i)}+p_{\sigma(i)}-1)\prod_{k=l+1}^r a^{(\sigma(k))}_n}{n^{w'_\sigma({\bf p},l)+q-|{\bf j}_\sigma|_l-t}}\nonumber\\
  &+\sum_{l=1}^{r} \sum_{\sigma\in \mathcal{S}_r}\sum_{|{\bf j}_\sigma|_l= w'_\sigma({\bf p},l)+q} (-1)^{|{\bf p}_\sigma|_l}C_l({\bf j}_\sigma;{\bf p}_\sigma)\nonumber\\
  &\quad\quad\quad\quad\quad\quad\quad\quad\quad\quad\times\prod_{i=1}^l D^{(A^{(\sigma(i))})}(j_{\sigma(i)}+p_{\sigma(i)}-1)\prod_{k=l+1}^r a^{(\sigma(k))}_0\nonumber\\
  &-2\sum_{l=0}^{r} \sum_{\sigma\in \mathcal{S}_r}\sum_{|{\bf j}_\sigma|_l+2t= w'_\sigma({\bf p},l)+q} (-1)^{|{\bf p}_\sigma|_l}C_l({\bf j}_\sigma;{\bf p}_\sigma)\nonumber\\
  &\quad\quad\quad\quad\quad\quad\quad\quad\quad\quad\times D\a(2t)\prod_{i=1}^l D^{(A^{(\sigma(i))})}(j_{\sigma(i)}+p_{\sigma(i)}-1)\prod_{k=l+1}^r a^{(\sigma(k))}_0\nonumber\\
  &=0,
\end{align}
where $w'_\sigma({\bf p},l):=|{\bf p}|-|{\bf p}_\sigma|_l+l$.
\end{thm}
\pf Applying the kernel $\pi\cot(\pi s;A)f(s;{\bf A})$ to the base function $r(s)=1/s^q$ with the help of identities (\ref{4.2}) and (\ref{4.3}), we arrive at the desired evaluation after
an elementary but rather tedious computation, which we leave to the interested reader.\hfill$\square$

It is clear that Theorems \ref{thm3.4} and \ref{thm3.6} are immediate corollaries of Theorems \ref{thm4.1} and \ref{thm4.2} with $r=2,p_1=m,p_2=p$ and $A^{(1)}=A,A^{(2)}=B$.

Further, letting all $A^{(l)}\in\{A_1,A_2\}\ (l=1,2,\ldots,r)$, we can get the following parity theorem.
\begin{thm}({\rm Parity theorem})\label{cor4.3}
(alternating) Euler sum $S_{e_1e_2\cdots e_r,g}\ (e_j\in\{p_j,{\bar p}_j\},g\in\{q,{\bar q}\})$ reduces to a combination of sums of lower orders whenever the weight $p_1+p_2+\cdots+p_r+q$ and the order $r$ are of the same parity.
\end{thm}

It should be emphasized that every (alternating) Euler sum of weight $w$ and degree $n$ is clearly a $\mathbb{Q}$-linear combination of (alternating) multiple zeta values of weight $w$ and depth less than or equal to $n+1$ (explicit formulas was given by Xu and Wang \cite{XW2018}). The multiple zeta values are defined by (\cite{H1992,DZ1994,Z2016})
\begin{align}\label{AMZV-Defn}
&\zeta( \mathbf{k})\equiv\zeta(k_1, \ldots, k_n):=\sum\limits_{m_1>\cdots>m_n\geq 1}\prod\limits_{j=1}^n m_j^{-|k_j|}{\rm sgn}(k_j)^{m_j},
\end{align}
where for convergence $|k_1|+\cdots+|k_j|> j$ for $j= 1, 2, \ldots, n$, and
\[{\rm sgn}(k_j):=\begin{cases}
   1  & \text{\;if\;} k_j>0,  \\
   -1, & \text{\;if\;} k_j<0.
\end{cases}\]
Here, we call $l(\mathbf{k}):=n$ and $|\mathbf{k}|:=\sum\limits_{j=1}^n|k_j|$ the depth and the weight of multiple zeta values, respectively. Further, we put a bar on top of
$-k_j$ if $k_j$ is negative integer. For example,
\begin{equation*}
\zeta(\bar 2,3,\bar 1,4)=\zeta(  -2,3,-1,4).
\end{equation*}
The systematic study of MZVs began in the early 1990s with the works of Hoffman \cite{H1992} and Zagier \cite{DZ1994}. Due to their surprising and sometimes mysterious appearance in the study of many branches of mathematics and theoretical physics, these special values have attracted a lot of attention and interest in the past three decades (for example, see the book by the
second author  \cite{Z2016}).
Rational relations among multiple zeta values are known through weight (sum of the indices) 22 and tabulated in the Multiple Zeta Value Data Mine (henceforth MZVDM) \cite{BBV2010}. Rational relations among alternating multiple zeta values are also tabulated in the MZVDM (through weight 12 \cite{BBV2010}).
With the help of results of \cite{BBV2010}, the author and Wang \cite{XW2018} have developed the Maple package to evaluate the non-alternating Euler sums of weight $2\leq w\leq 16$ and the alternating Euler sums of weight $1\leq w\leq 6$.

We record two examples to illustrate Theorem \ref{thm4.2} and Theorem \ref{cor4.3}.
\begin{exa} For integer $p>1$,
\begin{align*}
& \left( {{{\left( { - 1} \right)}^{p - 1}} - 1} \right){S_{{{\bar 1}^3},\bar p}} - 6{S_{\bar 1\bar 2,\bar p}} + 3\log (2)\left( {1 + {{\left( { - 1} \right)}^p}} \right){S_{{{\bar 1}^2},\bar p}} + 6\log (2){S_{\bar 2,\bar p}}\nonumber \\ & + 3\left( {\zeta \left( 2 \right) - \left( {1 + {{\left( { - 1} \right)}^p}} \right){{\log }^2}(2)} \right){S_{\bar 1,\bar p}}
 - 3{S_{\bar 3,\bar p}} + {\left( { - 1} \right)^p}{\log ^3}(2)\bar \zeta \left( p \right) - 3p{S_{{{\bar 1}^2},\overline {p + 1}}} - 3p{S_{\bar 2,\overline {p + 1}}} \nonumber
 \\&+ 6p\log (2){S_{\bar 1,\overline {p + 1}}} - p\left( {3{{\log }^2}(2) - \frac{1}{2}\zeta \left( 2 \right)} \right)\bar \zeta \left( {p + 1} \right)
  - 3\left( {\frac{{p\left( {p + 1} \right)}}{2} + {{\left( { - 1} \right)}^p}} \right){S_{\bar 1,\overline {p + 2}}}
  \nonumber \\&+ 3\left( {\frac{{p\left( {p + 1} \right)}}{2} + {{\left( { - 1} \right)}^p}} \right)\log (2)\bar \zeta \left( {p + 2} \right) + \left( {{{\log }^3}(2) - 3\log (2)\zeta \left( 2 \right) + \frac{9}{4}\zeta \left( 3 \right)} \right)\bar \zeta \left( p \right)\nonumber \\
 & + {\left( { - 1} \right)^p}\zeta \left( {p + 3} \right) - \frac{{p\left( {p + 1} \right)\left( {p + 2} \right)}}{6}\bar \zeta \left( {p + 3} \right)
 \nonumber \\&+ 3{\left( { - 1} \right)^p}\left( {{S_{{{\bar 1}^2},p + 1}} - 2\log (2){S_{\bar 1,p + 1}} + {{\log }^2}(2)\zeta \left( {p + 1} \right)} \right)\nonumber \\
 &+ {\rm{Res}}\left[ {{g_1}\left( s \right),s = 0} \right]=0,
\end{align*}
where
\begin{align*}
 {\rm{Res}}\left[ {{g_1}\left( s \right),s = 0} \right] =&  - 3\bar \zeta \left( {p + 3} \right) - 3\sum\limits_{\scriptstyle {j_1} + {j_2} = p + 1, \hfill \atop
  \scriptstyle {j_1},{j_2} \ge 0 \hfill} {\bar \zeta \left( {{j_1} + 1} \right)\bar \zeta \left( {{j_2} + 1} \right)} \nonumber \\
  &- \sum\limits_{\scriptstyle {j_1} + {j_2} + {j_3} = p, \hfill \atop
  \scriptstyle {j_1},{j_2},{j_3} \ge 0 \hfill} {\bar \zeta \left( {{j_1} + 1} \right)\bar \zeta \left( {{j_2} + 1} \right)\bar \zeta \left( {{j_3} + 1} \right)}  + \left( {1 - {{\left( { - 1} \right)}^p}} \right)\zeta \left( {p + 3} \right)\nonumber \\
  &+ 6\sum\limits_{\scriptstyle {j_1} + 2{j_2} = p + 2, \hfill \atop
  \scriptstyle {j_1} \ge 0,{j_2} \ge 1 \hfill} {\bar \zeta \left( {{j_1} + 1} \right)\zeta \left( {2{j_2}} \right)} \nonumber\\& + 6\sum\limits_{\scriptstyle {j_1} + {j_2} + 2{j_3} = p + 1, \hfill \atop
  \scriptstyle {j_1},{j_2} \ge 0,{j_3} \ge 1 \hfill} {\bar \zeta \left( {{j_1} + 1} \right)\bar \zeta \left( {{j_2} + 1} \right)\zeta \left( {2{j_3}} \right)} \nonumber \\
 & + 2\sum\limits_{\scriptstyle {j_1} + {j_2} + {j_3} + 2{j_4} = p, \hfill \atop
  \scriptstyle {j_1},{j_2},{j_3} \ge 0,{j_4} \ge 1 \hfill} {\bar \zeta \left( {{j_1} + 1} \right)\bar \zeta \left( {{j_2} + 1} \right)\bar \zeta \left( {{j_3} + 1} \right)\zeta \left( {2{j_4}} \right)}.
\end{align*}
\end{exa}
\begin{exa} For integer $p>1$,
\begin{align*}
& \left( {{{\left( { - 1} \right)}^{p - 1}} - 1} \right){S_{{{\bar 1}^3},p}}{\rm{ + }}6{S_{\bar 1\bar 2,p}} - 3\log (2)\left( {1 + {{\left( { - 1} \right)}^p}} \right){S_{{{\bar 1}^2},p}} - 6\log (2){S_{\bar 2,p}}\nonumber \\
 & + 3\left( {\left( {1 + {{\left( { - 1} \right)}^p}} \right){{\log }^2}(2)+2 \zeta \left( 2 \right)} \right){S_{\bar 1,p}} + 3{S_{\bar 3,p}} - {\left( { - 1} \right)^p}{\log ^3}(2)\zeta \left( p \right) + 3p{S_{{{\bar 1}^2},p + 1}} + 3p{S_{\bar 2,p + 1}}\nonumber \\
&  - 6p\log (2){S_{\bar 1,p + 1}} + p\left( {3{{\log }^2}(2) - \frac{5}{2}\zeta \left( 2 \right)} \right)\zeta \left( {p + 1} \right) + 3\left( {\frac{{p\left( {p + 1} \right)}}{2} + {{\left( { - 1} \right)}^p}} \right){S_{\bar 1,p + 2}} \nonumber  \\
&- 3\left( {\frac{{p\left( {p + 1} \right)}}{2} + {{\left( { - 1} \right)}^p}} \right)\log (2)\zeta \left( {p + 2} \right) - \left( {{{\log }^3}(2) + 6\log (2)\zeta \left( 2 \right) + \frac{9}{4}\zeta \left( 3 \right)} \right)\zeta \left( p \right)\nonumber \\
 & - {\left( { - 1} \right)^p}\bar \zeta \left( {p + 3} \right) + \frac{{p\left( {p + 1} \right)\left( {p + 2} \right)}}{6}\zeta \left( {p + 3} \right)\nonumber \\
 &- 3{\left( { - 1} \right)^p}\left( {{S_{{{\bar 1}^2},\overline{p + 1}}} - 2\log (2){S_{\bar 1,\overline{p + 1}}} + {{\log }^2}(2)\bar \zeta \left( {p + 1} \right)} \right)\nonumber \\
 & + {\rm{Res}}\left[ {{g_2}\left( s \right),s = 0} \right] = 0 ,
\end{align*}
where
\begin{align*}
 {\rm{Res}}\left[ {{g_2}\left( s \right),s = 0} \right] = & - 3\bar \zeta \left( {p + 3} \right) - 3\sum\limits_{\scriptstyle {j_1} + {j_2} = p + 1, \hfill \atop
  \scriptstyle {j_1},{j_2} \ge 0 \hfill} {\bar \zeta \left( {{j_1} + 1} \right)\bar \zeta \left( {{j_2} + 1} \right)}\nonumber  \\
 & - \sum\limits_{\scriptstyle {j_1} + {j_2} + {j_3} = p, \hfill \atop
  \scriptstyle {j_1},{j_2},{j_3} \ge 0 \hfill} {\bar \zeta \left( {{j_1} + 1} \right)\bar \zeta \left( {{j_2} + 1} \right)\bar \zeta \left( {{j_3} + 1} \right)}  - \left( {1 - {{\left( { - 1} \right)}^p}} \right)\bar \zeta \left( {p + 3} \right) \nonumber \\
 & - 6\sum\limits_{\scriptstyle {j_1} + 2{j_2} = p + 2, \hfill \atop
  \scriptstyle {j_1} \ge 0,{j_2} \ge 1 \hfill} {\bar \zeta \left( {{j_1} + 1} \right)\bar \zeta \left( {2{j_2}} \right)}\nonumber   \\
 & - 6\sum\limits_{\scriptstyle {j_1} + {j_2} + 2{j_3} = p + 1, \hfill \atop
  \scriptstyle {j_1},{j_2} \ge 0,{j_3} \ge 1 \hfill} {\bar \zeta \left( {{j_1} + 1} \right)\bar \zeta \left( {{j_2} + 1} \right)\bar \zeta \left( {2{j_3}} \right)} \nonumber  \\
  &- 2\sum\limits_{\scriptstyle {j_1} + {j_2} + {j_3} + 2{j_4} = p, \hfill \atop
  \scriptstyle {j_1},{j_2},{j_3} \ge 0,{j_4} \ge 1 \hfill} {\bar \zeta \left( {{j_1} + 1} \right)\bar \zeta \left( {{j_2} + 1} \right)\bar \zeta \left( {{j_3} + 1} \right)\bar \zeta \left( {2{j_4}} \right)} .
\end{align*}
\end{exa}

Applying Au's Mathematica package \cite{Au2020}, we can get the following specific cases.
\begin{align*}
S_{\bar1 1 \bar 2,2}&=-\frac{9}{4} \z(\bar 5,1)+\frac{1}{6} \pi ^2 \text{Li}_4\left(\frac{1}{2}\right)-4 \text{Li}_6\left(\frac{1}{2}\right)+2 \text{Li}_4\left(\frac{1}{2}\right) \log ^2(2)+\frac{35 \zeta (3)^2}{32}\\&\quad+\frac{7}{6} \zeta (3) \log ^3(2)-\frac{1}{3} \pi ^2 \zeta (3) \log (2)+\frac{279}{64} \zeta (5) \log (2)+\frac{187 \pi ^6}{362880}\\&\quad+\frac{7 \log ^6(2)}{90}-\frac{1}{16} \pi ^2 \log ^4(2)-\frac{1}{144} \pi ^4 \log ^2(2),\\
S_{\bar1 1 2,2}&=-\frac{13}{4}\z(\bar 5,1)-\frac{1}{12} \pi ^2 \text{Li}_4\left(\frac{1}{2}\right)+12 \text{Li}_6\left(\frac{1}{2}\right)+2 \text{Li}_4\left(\frac{1}{2}\right) \log ^2(2)\\&\quad+8 \text{Li}_5\left(\frac{1}{2}\right) \log (2)+\frac{45 \zeta (3)^2}{32}+\frac{1}{4} \pi ^2 \zeta (3) \log (2)+\frac{31}{8} \zeta (5) \log (2)\\&\quad-\frac{319 \pi ^6}{22680}+\frac{\log ^6(2)}{30}-\frac{5}{288} \pi ^2 \log ^4(2)+\frac{1}{288} \pi ^4 \log ^2(2),\\
S_{\bar1 \bar1 \bar2,2}&=16 \z(\bar 5,1)-\frac{1}{6} \pi ^2 \text{Li}_4\left(\frac{1}{2}\right)-32 \text{Li}_6\left(\frac{1}{2}\right)-8 \text{Li}_5\left(\frac{1}{2}\right) \log (2)-\frac{43 \zeta (3)^2}{8}\\&\quad-\frac{7}{6} \zeta (3) \log ^3(2)-\frac{3}{8} \pi ^2 \zeta (3) \log (2)-\frac{93}{4} \zeta (5) \log (2)+\frac{15833 \pi ^6}{362880}\\&\quad+\frac{\log ^6(2)}{45}-\frac{1}{144} \pi ^2 \log ^4(2)+\frac{103}{720} \pi ^4 \log ^2(2),\\\
S_{\bar1 \bar1 \bar2,\bar2}&=12 \z(\bar 5,1)-24 \text{Li}_6\left(\frac{1}{2}\right)+4 \text{Li}_4\left(\frac{1}{2}\right) \log ^2(2)-\frac{129 \zeta (3)^2}{32}-\frac{11}{48} \pi ^2 \zeta (3) \log (2)\\&\quad-\frac{93}{4} \zeta (5) \log (2)+\frac{5669 \pi ^6}{181440}+\frac{2 \log ^6(2)}{15}-\frac{1}{12} \pi ^2 \log ^4(2)+\frac{49}{360} \pi ^4 \log ^2(2),\\
S_{\bar1\bar1\bar1,3}&=3 \z(\bar 5,1)+12 \text{Li}_5\left(\frac{1}{2}\right) \log (2)-\frac{27 \zeta (3)^2}{16}+\frac{7}{4} \zeta (3) \log ^3(2)+\frac{9}{16} \pi ^2 \zeta (3) \log (2)\\&\quad-\frac{279}{16} \zeta (5) \log (2)+\frac{97 \pi ^6}{40320}-\frac{1}{10} \log ^6(2)+\frac{1}{6} \pi ^2 \log ^4(2)+\frac{19}{240} \pi ^4 \log ^2(2),\\
S_{\bar1\bar1\bar1,\bar3}&=9\z(\bar 5,1)-12 \text{Li}_6\left(\frac{1}{2}\right)-6 \text{Li}_4\left(\frac{1}{2}\right) \log ^2(2)-\frac{285 \zeta (3)^2}{64}+\frac{9}{16} \pi ^2 \zeta (3) \log (2)\\&\quad-\frac{279}{16} \zeta (5) \log (2)+\frac{2 \pi ^6}{105}-\frac{1}{15} 4 \log ^6(2)+\frac{7}{24} \pi ^2 \log ^4(2)+\frac{19}{240} \pi ^4 \log ^2(2),\\
S_{\bar1\bar1 1,3}&=\frac{11}{2} \z(\bar 5,1)-8 \text{Li}_6\left(\frac{1}{2}\right)-\frac{3 \zeta (3)^2}{2}-\frac{7}{6} \zeta (3) \log ^3(2)+\frac{1}{12} \pi ^2 \zeta (3) \log (2)\\&\quad-\frac{341}{32} \zeta (5) \log (2)+\frac{131 \pi ^6}{12096}-\frac{1}{90} \log ^6(2)+\frac{1}{36} \pi ^2 \log ^4(2)+\frac{4}{45} \pi ^4 \log ^2(2),\\
S_{\bar1\bar1 1,\bar3}&=4 \z(\bar 5,1)-\frac{1}{6} \pi ^2 \text{Li}_4\left(\frac{1}{2}\right)-4 \text{Li}_6\left(\frac{1}{2}\right)+2 \text{Li}_4\left(\frac{1}{2}\right) \log ^2(2)+4 \text{Li}_5\left(\frac{1}{2}\right) \log (2)\\&\quad-\frac{83 \zeta (3)^2}{64}-\frac{7}{12} \zeta (3) \log ^3(2)-\frac{1}{16} \pi ^2 \zeta (3) \log (2)-\frac{341}{32} \zeta (5) \log (2)+\frac{121 \pi ^6}{15120}\\&\quad+\frac{2 \log ^6(2)}{45}-\frac{1}{48} \pi ^2 \log ^4(2)+\frac{23}{240} \pi ^4 \log ^2(2),\\
S_{111,\bar3}&=\frac{9}{2} \z(\bar 5,1)+\frac{1}{2} \pi ^2 \text{Li}_4\left(\frac{1}{2}\right)-12 \text{Li}_6\left(\frac{1}{2}\right)-6 \text{Li}_4\left(\frac{1}{2}\right) \log ^2(2)-12 \text{Li}_5\left(\frac{1}{2}\right) \log (2)\\&\quad-\frac{207 \zeta (3)^2}{64}-\frac{7}{4} \zeta (3) \log ^3(2)+\frac{7}{16} \pi ^2 \zeta (3) \log (2)+\frac{257 \pi ^6}{20160}\\&\quad-\frac{1}{6} \log ^6(2)+\frac{7}{48} \pi ^2 \log ^4(2)-\frac{1}{48} \pi ^4 \log ^2(2).
\end{align*}

\section{Parameterized Extensions of Ramanujan-Type Identities}

Infinite series involving the trigonometric and hyperbolic functions have attracted the attention of many authors. Ramanujan evaluate many infinite series of
hyperbolic functions in his notebooks \cite{R2012} and lost notebook \cite{R1988}. Andrews and Berndt's books \cite{AB2009,AB2013} contain many such results well as numerous references. In this section, we provide some parameterized extensions of Ramanujan-type identities
that involve hyperbolic series.

We define the hyperbolic cotangent function with sequence $A$ by
\begin{align}
\coth \left( {s;A} \right): = i\cot \left( {is;A} \right).
\end{align}

\begin{thm}\label{thm5.1} Let $k>1$ be a integer and $\alpha,\beta$ be reals with $\alpha,\beta\neq 0$ with $A,B$ defined above, then
\begin{align}\label{5.1}
 &{\alpha ^{2k - 1}}\beta \pi \sum\limits_{n = 1}^\infty  {\frac{{{a_n}\coth \left( {\frac{{n\beta \pi }}{\alpha };B} \right)}}{{{n^{2k - 1}}}}}  + {\left( { - 1} \right)^k}{\beta ^{2k - 1}}\alpha \pi \sum\limits_{n = 1}^\infty  {\frac{{{b_n}\coth \left( {\frac{{n\alpha \pi }}{\beta };A} \right)}}{{{n^{2k - 1}}}}} \nonumber \\
  &= {b_0}{\alpha ^{2k}}D\a(2k)  + {\left( { - 1} \right)^k}{a_0}{\beta ^{2k}}D\B(2k)  - 2\sum\limits_{\scriptstyle {j_1} + {j_2} = k, \hfill \atop
  \scriptstyle {j_1},{j_2} \ge 1 \hfill} {{{\left( { - 1} \right)}^{{j_2}}}{\alpha ^{2{j_1}}}{\beta ^{2{j_2}}} D\a(2j_1)D\B(2j_2)}.
\end{align}
\end{thm}
\pf To prove (\ref{5.1}), we consider the function
\[f\left( {s;A,B} \right): = \frac{{\pi^2 \cot \left( {\pi \alpha s;A} \right)\coth \left( {\pi \beta s;B} \right)}}{{{s^{2k-1}}}}.\]
It is obvious that the function $f\left( {s;A,B} \right)$ is meromorphic in the entire complex plane with a simple pole at $s=0$ and simple poles $s=\frac{n}{\alpha }$ and $s=\frac{n}{\beta }i$ for each integer $n$. Then, brief calculations show that
\begin{align}
&\begin{aligned}\label{5.2}
 {\rm{Res}}\left[ {f\left( {s;A,B} \right),s = \frac{n}{\alpha }} \right] = {\alpha ^{2k - 2}}\pi \frac{{{a_n}\coth \left( {\frac{{n\beta \pi }}{\alpha };B} \right)}}{{{n^{2k - 1}}}},
\end{aligned}\\
&\begin{aligned}\label{5.3}
 {\rm{Res}}\left[ {f\left( {s;A,B} \right),s = \frac{n}{\beta }i} \right] = {\left( { - 1} \right)^k}{\beta ^{2k - 2}}\pi \frac{{{b_n}\coth \left( {\frac{{n\alpha \pi }}{\beta };A} \right)}}{{{n^{2k - 1}}}},
\end{aligned}\\
&\begin{aligned}\label{5.4}
 {\rm{Res}}\left[ {f\left( {s;A,B} \right),s = 0} \right] =  -2{b_0}\frac{\alpha ^{2k-1}}{\beta}D\a(2k)  -2{\left( { - 1} \right)^k}{a_0}\frac{\beta ^{2k-1}}{\alpha}D\B(2k) \\ +4\sum\limits_{\scriptstyle {j_1} + {j_2} = k, \hfill \atop
  \scriptstyle {j_1},{j_2} \ge 1 \hfill} {{{\left( { - 1} \right)}^{{j_2}}}{\alpha ^{2{j_1}-1}}{\beta ^{2{j_2}-1}} D\a(2j_1)D\B(2j_2)}
\end{aligned}
\end{align}
for each integer $n\ (n\neq 0)$. Hence, using (\ref{5.2})-(\ref{5.4}) and the residue theorem, we reach the desired result. \hfill$\square$

An elementary calculation, the ({\ref{5.1}}) can be rewritten in the form
\begin{align}\label{5.5}
& \alpha {\beta ^k}\sum\limits_{n = 1}^\infty  {\frac{{{a_n}\coth \left( {n\alpha ;B} \right)}}{{{n^{2k - 1}}}}}  + {\left( { - 1} \right)^k}\beta {\alpha ^k}\pi \sum\limits_{n = 1}^\infty  {\frac{{{b_n}\coth \left( {n\beta ;A} \right)}}{{{n^{2k - 1}}}}} \nonumber \\
 & = {b_0}{\beta ^k}D\a(2k)  + {\left( { - 1} \right)^k}{a_0}{\alpha ^k}D\B(2k) - 2\sum\limits_{\scriptstyle {j_1} + {j_2} = k, \hfill \atop
  \scriptstyle {j_1},{j_2} \ge 1 \hfill} {{{\left( { - 1} \right)}^{{j_2}}}{\alpha ^{{j_2}}}{\beta ^{{j_1}}}D\a(2j_1)D\B(2j_2)},
\end{align}
where $\alpha$ and $\beta$ are real numbers with $\alpha \beta  = {\pi ^2}$.

If setting $A,B\in\{A_1,A_2\}$, then we get the following well-known results.
\begin{cor}\label{cor5.3} If $k>1$ is a positive integer number and $\alpha,\beta$ are real numbers such that $\alpha \beta  = {\pi ^2}$, then
\begin{align}\label{5.6}
&\alpha {\beta ^k}\sum\limits_{n = 1}^\infty  {\frac{{\coth \left( {n\alpha } \right)}}
{{{n^{2k - 1}}}}}  + {\left( { - 1} \right)^k}{\alpha ^k}\beta \sum\limits_{n = 1}^\infty  {\frac{{\coth \left( {n\beta } \right)}}
{{{n^{2k - 1}}}}}\nonumber \\
& = \left( {{\beta ^k} + \left( { - 1} \right)^k{\alpha ^k}} \right)\zeta \left( {2k} \right) - 2\sum\limits_{\scriptstyle {k_1} + {k_2} = k,\hfill \atop \scriptstyle {k_1}, {k_2} \geqslant 1\hfill}  {{{\left( { - 1} \right)}^{{k_2}}}{\alpha ^{{k_2}}}{\beta ^{{k_1}}}\zeta \left( {2{k_1}} \right)\zeta \left( {2{k_2}} \right)}.
\end{align}
\end{cor}
\begin{cor} \label{cor5.4} If $k$ is a positive integer number and $\alpha,\beta$ are real numbers such that $\alpha \beta  = {\pi ^2}$, then
\begin{align}\label{5.7}
&\alpha {\beta ^k}\sum\limits_{n = 1}^\infty  {\frac{{{{\left( { - 1} \right)}^{n - 1}}}}
{{{n^{2k - 1}}\sinh \left( {n\alpha } \right)}}}  + {\left( { - 1} \right)^k}{\alpha ^k}\beta \sum\limits_{n = 1}^\infty  {\frac{{{{\left( { - 1} \right)}^{n - 1}}}}
{{{n^{2k - 1}}\sinh \left( {n\beta } \right)}}}\nonumber \\
& = \left( {{\beta ^k} + \left( { - 1} \right)^k{\alpha ^k}} \right)\bar \zeta \left( {2k} \right) + 2\sum\limits_{\scriptstyle {k_1} + {k_2} = k,\hfill \atop \scriptstyle {k_1}, {k_2} \geqslant 1\hfill}{{{\left( { - 1} \right)}^{{k_2}}}{\alpha ^{{k_2}}}{\beta ^{{k_1}}}\bar \zeta \left( {2{k_1}} \right)\bar \zeta \left( {2{k_2}} \right)}.
\end{align}
\end{cor}
\begin{cor}\label{cor5.5} If $k$ is a positive integer number and $\alpha,\beta$ are real numbers such that $\alpha \beta  = {\pi ^2}$, then
\begin{align}\label{5.8}
&\alpha {\beta ^k}\sum\limits_{n = 1}^\infty  {\frac{1}
{{{n^{2k - 1}}\sinh \left( {n\alpha } \right)}}}  + {\left( { - 1} \right)^k}{\alpha ^k}\beta \sum\limits_{n = 1}^\infty  {\frac{{\coth \left( {n\beta } \right)}}
{{{n^{2k - 1}}}}{{\left( { - 1} \right)}^n}}\nonumber \\
&= \left( {{\beta ^k} - {{\left( { - 1} \right)}^k}\left( {1 - {2^{1 - 2k}}} \right){\alpha ^k}} \right)\zeta \left( {2k} \right) + 2\sum\limits_{\scriptstyle {k_1} + {k_2} = k,\hfill \atop \scriptstyle {k_1}, {k_2} \geqslant 1\hfill} {{{\left( { - 1} \right)}^{{k_2}}}{\alpha ^{{k_2}}}{\beta ^{{k_1}}}\zeta \left( {2{k_1}} \right)\bar \zeta \left( {2{k_2}} \right)} .
\end{align}
\end{cor}

Applying the fact
\[\coth x = 1 + \frac{2}
{{{e^{2x}} - 1}}\]
to (\ref{5.6}),
we obtain the well-known Ramanujan's formula for $\zeta(2k-1)\ (k>1)$
\begin{align}\label{5.9}
&{\left( {4\beta } \right)^{ - \left( {k - 1} \right)}}\left\{ {\frac{1}
{2}\zeta \left( {2k - 1} \right) + \sum\limits_{n = 1}^\infty  {\frac{1}
{{{n^{2k - 1}}\left( {{e^{2n\alpha }} - 1} \right)}}} } \right\}\nonumber\\
& \quad- {\left( { - 4\alpha } \right)^{ - \left( {k - 1} \right)}}\left\{ {\frac{1}
{2}\zeta \left( {2k - 1} \right) + \sum\limits_{n = 1}^\infty  {\frac{1}
{{{n^{2k - 1}}\left( {{e^{2n\beta }} - 1} \right)}}} } \right\}\nonumber\\
& = \sum\limits_{j = 0}^k {{{\left( { - 1} \right)}^{j - 1}}\frac{{{B_{2j}}{B_{2k - 2j}}}}
{{\left( {2j} \right){\text{!}}\left( {2k - 2j} \right)!}}{\alpha ^j}{\beta ^{k - j}}}.
\end{align}
The Ramanujan's formula for $\zeta(2k-1)$ (\ref{5.9}) appears as Entry 21(i) in Chapter 14 of Ramanujan's second notebook \cite{R2012}. It also appears in a formerly unpublished manuscript of Ramanujan that was published in its original handwritten form with his lost notebook \cite{R1988}. The first published proof of (\ref{5.9}) is due to Malurkar \cite{M1925} in 1925-1926. This partial manuscript was initially examined in detail by the Berndt in \cite{B2004}, and by Andrews and Berndt in their fourth book on Ramanujan's lost notebook \cite{AB2013}.
There are many results in the literature generalizing or extending the result. The two most extensive papers in this direction are perhaps those by Bradley \cite{Br2002} and Sitaramachandrarao \cite{S1}.

\begin{thm} Let $\alpha,\beta,\gamma$ be reals with $\alpha,\beta,\gamma\neq 0$, we have
\begin{align}\label{5.10}
&\alpha^{2k}\beta\gamma\pi^2\sum\limits_{n=1}^\infty \frac{a_n\coth(\frac{\pi \beta n}{\alpha};B)\coth(\frac{\pi \gamma n}{\alpha};C)}{n^{2k}}\nonumber\\
&+(-1)^{k-1}\beta^{2k}\alpha\gamma \pi^2\sum\limits_{n=1}^\infty \frac{b_n\coth(\frac{\pi \alpha n}{\beta};A)\cot(\frac{\pi \gamma n}{\beta};C)}{n^{2k}}\nonumber\\
&+(-1)^{k-1}\gamma^{2k}\alpha\beta\pi^2\sum\limits_{n=1}^\infty \frac{c_n\coth(\frac{\pi \alpha n}{\gamma};A)\cot(\frac{\pi \beta n}{\gamma};B)}{n^{2k}}\nonumber\\
&-2a_0(-1)^k\sum\limits_{\scriptstyle j_1+j_2=k+1, \hfill \atop
  \scriptstyle j_1,j_2\geq 1 \hfill}D\B(2j_1)D\C(2j_2)\beta^{2j_1}\gamma^{2j_2}\nonumber\\
&+2b_0\sum\limits_{\scriptstyle j_1+j_2=k+1, \hfill \atop
  \scriptstyle j_1,j_2\geq 1 \hfill} (-1)^{j_2}D\a(2j_1)D\C(2j_2)\alpha^{2j_1}\gamma^{2j_2}\nonumber\\
&+2c_0\sum\limits_{\scriptstyle j_1+j_2=k+1, \hfill \atop
  \scriptstyle j_1,j_2\geq 1 \hfill} (-1)^{j_2}D\a(2j_1)D\B(2j_2)\alpha^{2j_1}\beta^{2j_2}\nonumber\\
&-4\sum\limits_{\scriptstyle j_1+j_2+j_3=k+1, \hfill \atop
  \scriptstyle j_1,j_2,j_3\geq 1 \hfill} (-1)^{j_2+j_3}D\a(2j_1)D\B(2j_2)D\C(2j_3)\alpha^{2j_1}\beta^{2j_2}\gamma^{2j_3}\nonumber\\
&+a_0b_0\gamma^{2k+2}(-1)^kD\C(2k+2)+a_0c_0\beta^{2k+2}(-1)^kD\B(2k+2)\nonumber\\&-b_0c_0\alpha^{2k+2}D\a(2k+2)=0.
\end{align}
\end{thm}
\pf Consider
\[f(s;A,B,C):=\frac{\pi^3\cot(\pi\alpha s;A)\coth(\pi \beta s;B)\coth(\pi \gamma s;C)}{s^{2k}}\quad (s\in \mathbb{C}\ {\rm and}\ \frac{\gamma}{\beta} \notin {\mathbb{Q}}).\]
Then, proceeding in the same fashion as in the proof of Theorem \ref{thm5.1}, we obtain the desired formula. \hfill$\square$

 Setting $a_n=c_n=1,b_n=(-1)^n$ and $k=\alpha=1,\beta=\sqrt 2,\gamma=2$ in (\ref{5.10}) yields
\[\sum\limits_{n = 1}^\infty  {\frac{{\coth \left( {2\pi n} \right)}}{{{n^2}\sinh \left( {\sqrt 2 \pi n} \right)}}}  + \sqrt 2 \sum\limits_{n = 1}^\infty  {\frac{{\coth \left( {\frac{{\pi n}}{{\sqrt 2 }}} \right)\cot \left( {\sqrt 2 \pi n} \right)}}{{{n^2}}}{{\left( { - 1} \right)}^n}}  + 2\sum\limits_{n = 1}^\infty  {\frac{{\coth \left( {\frac{{\pi n}}{2}} \right)}}{{{n^2}\sin \left( {\frac{{\pi n}}{{\sqrt 2 }}} \right)}}}  = \frac{{97\sqrt 2 }}{{120}}\zeta \left( 2 \right).\]
In general, let $A:=\{a_n\}$ and $B_l:=\{b^{(l)}_n\}\ (l=1,2,\ldots,m;m\in\N)$  be sequences of complex numbers (if $n\rightarrow \infty$, then ${a_n} = o\left( {{n^{s_1} }}\right)$ and $ b^{(l)}_n=o\left(n^{t_l}\right)$ with $s_1,t_l<1$), $k$ be a positive integer and $\alpha,\beta_l\ (l=1,2,\ldots,m;m\in\N)$ be reals with $\alpha,\beta_l\neq 0$ and $\frac {\beta_l}{\beta_j} \notin {\mathbb{Q}}\ (l\neq j)$, define
\[f\left( {s;A,{B_1},{B_2}, \cdots ,{B_m}} \right): = \frac{{{\pi ^{m + 1}}\cot \left( {\pi \alpha s;A} \right)\prod\limits_{l = 1}^m {\coth \left( {\pi {\beta _l}s;{B_l}} \right)} }}{{{s^{2k + \left[ {1 - {{\left( { - 1} \right)}^m}} \right]/2}}}}\quad (s\in \mathbb{C}),\]
by elementary residue calculations, we can get the following general conclusion:
\begin{align}\label{5.11}
&{\left( { - 1} \right)^{k + m/2 - \left[ {1 - {{\left( { - 1} \right)}^m}} \right]/4}}{\pi ^m}\sum\limits_{j = 1}^m {{\beta_j^{2k - \left[ {3 - {{\left( { - 1} \right)}^m}} \right]/2}}\sum\limits_{n = 1}^\infty  {\frac{{b_n^{\left( j \right)}\coth \left( {\frac{{\pi \alpha n}}{{{\beta _j}}};A} \right)\prod\limits_{l = 1,l \ne j}^m {\cot \left( {\frac{{\pi {\beta _l}n}}{{{\beta _j}}};{B_l}} \right)} }}{{{n^{2k - \left[ {1 - {{\left( { - 1} \right)}^m}} \right]/2}}}}} }  \nonumber\\
& + {\alpha ^{2k - \left[ {3 - {{\left( { - 1} \right)}^m}} \right]/2}}{\pi ^m}\sum\limits_{n = 1}^\infty  {\frac{{{a_n}\prod\limits_{l = 1}^m {\coth \left( {\frac{{\pi {\beta _l}n}}{\alpha };{B_l}} \right)} }}{{{n^{2k - \left[ {1 - {{\left( { - 1} \right)}^m}} \right]/2}}}}} + \frac{1}{2}{\rm Res}[f\left( {s;A,{B_1}, \cdots ,{B_m}} \right),s=0]=0.
\end{align}

Next, we give a partial fraction decomposition for $\frac{\pi}{2}\cot(\sqrt{w\alpha};A)\coth(\sqrt{w\beta};B)$. From the definition of $\coth \left( {s;A} \right)$, using (\ref{2.6}), we have
\begin{align}
\pi \coth(\pi s; A)=\frac{a_{|n|}}{s-ni}-\sj (-i\sigma_{-n})^j R\a_{|n|}(j)(s-ni)^{j-1},
\end{align}
where $|s-ni|<1$ with $s\neq ni$ and $n$ is any integer. If $n=0$, then for $s\rightarrow 0$,
\begin{align}\label{5.14}
\pi \coth(\pi s; A)=\frac{a_0}{s}-2\sj (-1)^j D\a(2j)s^{2j-1}.
\end{align}

Now, we use the residue computations to establish a Sitaramachandrarao-type formula for $\pi^2xy \cot(\pi x;A)\coth(\pi y;B)$.
\begin{thm}\label{thm5.6} Let $x$ and $y$ be reals with $x,y\neq 0$, then
\begin{align}\label{5.15}
&\pi^2xy \cot(\pi x;A)\coth(\pi y;B)\nonumber\\&=a_0b_0+2a_0D\B(2)y^2-2b_0D\a(2)x^2\nonumber\\
&\quad-2\pi xy \su \left\{\frac{x^2a_n \coth(\pi ny/x;B)}{n(n^2-x^2)}+\frac{y^2b_n \coth(\pi nx/y;A)}{n(n^2+y^2)}\right\}.
\end{align}
\end{thm}
\pf Consider
\[F(s)=\frac{\pi^2\cot(\pi xs;A)\coth(\pi ys;B) }{s(s^2-1)},\]
which has simple poles at $s=\pm\frac{n}{x},\pm \frac{n}{y}i$ and $\pm 1$\ $(n\in \N)$, with residues
\begin{align*}
&{\rm Res}[F(s),s=\pm n/x]=\pi\frac{x^2a_n \coth(\pi ny/x;B)}{n(n^2-x^2)},\\
&{\rm Res}[F(s),s=\pm ni/y]=\pi\frac{y^2b_n \coth(\pi nx/y;A)}{n(n^2+y^2)},\\
&{\rm Res}[F(s),s=\pm 1]=\pi^2\frac{\cot(\pi x;A)\coth(\pi y;B)}{2}.
\end{align*}
Clearly $F(s)$ also has a triple pole at $s=0$. Using related formulas, if $s\rightarrow 0$, then we have
\begin{align*}
F(s)=\frac{1}{s^3} \frac{1}{xy(s^2-1)} \left\{a_0b_0+2a_0D\B(2)(ys)^2-2b_0D\a(2)(xs)^2+o(s^2)\right\}.
\end{align*}
Hence, we can find the residue
\begin{align*}
{\rm Res}[F(s),s=0]=-\frac{a_0b_0}{xy}-2a_0D\B(2)\frac{y}{x}+2b_0D\a(2)\frac{x}{y}.
\end{align*}
Finally, summing these four contributions yields the statement of the theorem.
\hfill$\square$

Obviously, if $A=B=A_1$ in Theorem \ref{thm5.6}, we obtain the Sitaramachandrarao's formula
\begin{align}\label{5.16}
{\pi ^2}xy\cot \left( {\pi x} \right)\coth \left( {\pi y} \right) =& 1 + 2\zeta( 2)\left( {{y^2} - {x^2}} \right)\nonumber \\&- 2\pi xy\sum\limits_{n = 1}^\infty  {\left\{ {\frac{{{x^2}\coth \left( {\frac{{\pi ny}}
{x}} \right)}}
{{n\left( {{n^2} - {x^2}} \right)}} + \frac{{{y^2}\coth \left( {\frac{{\pi nx}}
{y}} \right)}}
{{n\left( {{n^2} - {y^2}} \right)}}} \right\}},
\end{align}
which can also be found in  Sitaramachandrarao\cite{S1} and Berndt \cite[pp. 271-272]{B1989}. Hence, (\ref{5.15}) can be regarded as a parameterized version of Sitaramachandrarao's formula (\ref{5.16}).
\begin{thm} If $\alpha$ and $\beta$ are positive numbers such that $\alpha \beta  = {\pi ^2}$, and $w>0$, then
\begin{align}\label{5.17}
&\frac{\pi}{2}\cot(\sqrt{w\alpha};A)\coth(\sqrt{w\beta};B)\nonumber\\
&=\frac{a_0b_0}{2w}+\frac{a_0\beta D\B(2)-b_0\alpha D\a(2)}{\pi^2}\nonumber\\
&\quad+\su \left\{\frac{n\alpha\coth(n\alpha;A)}{w+n^2\alpha}b_n+\frac{n\beta\coth(n\beta;B)}{w-n^2\beta}a_n \right\}\nonumber\\
&\quad+\su \frac{a_n\coth(n\beta;B)-b_n\coth(n\alpha;A)}{n}.
\end{align}
\end{thm}
\pf Let $\pi x=\sqrt{w\alpha}$ and $\pi y=\sqrt{w\beta}$ in (\ref{5.15}) to deduce that
\begin{align*}
&\frac{\pi}{2}\cot(\sqrt{w\alpha};A)\coth(\sqrt{w\beta};B)\nonumber\\
&=\frac{a_0b_0}{2w}+\frac{a_0\beta D\B(2)-b_0\alpha D\a(2)}{\pi^2}\nonumber\\
&\quad-w\su \left\{\frac{a_n\coth(n\beta;B)}{n(n^2\beta-w)}+\frac{b_n\coth(n\alpha;A)}{n(n^2\alpha+w)}\right\}.
\end{align*}
Then using the elementary identities
\begin{align*}
\frac{1}{n(n^2\beta-w)}=\frac{1}{w}\left\{\frac{n\beta}{n^2\beta-w}-\frac{1}{n}\right\}
\end{align*}
and
\begin{align*}
\frac{1}{n(n^2\alpha+w)}=\frac{1}{w}\left\{\frac{1}{n}-\frac{n\alpha}{n^2\alpha+w}\right\},
\end{align*}
With a direct calculation, we complete this proof. \hfill$\square$

Finally, we end this paper by two corollaries.

\begin{cor}(\cite{AB2013}) If $\alpha$ and $\beta$ are positive numbers such that $\alpha \beta  = {\pi ^2}$, and $w>0$, then
\begin{align}\label{5.18}
&\frac{\pi}{2}\cot(\sqrt{w\alpha})\coth(\sqrt{w\beta})=\frac{1}{2w}+\frac{1}{2}\log\left(\frac{\beta}{\alpha}\right)+\su \left\{\frac{n\alpha\coth(n\alpha)}{w+n^2\alpha}+\frac{n\beta\coth(n\beta)}{w-n^2\beta} \right\}.
\end{align}
\end{cor}
\pf Setting $A=B=A_1$ in (\ref{5.17}) and using the identity \cite[pp.274]{AB2013}
\[\su \frac{1}{n(e^{2n\alpha}-1)}-\frac{1}{4}\log\alpha +\frac{\alpha}{12}=\su \frac{1}{n(e^{2n\beta}-1)}-\frac{1}{4}\log\beta +\frac{\beta}{12},\]
and noting that the fact
\[\coth(x)=1+\frac{2}{e^{2x}-1}\]
we complete the proof of (\ref{5.18}).\hfill$\square$

\begin{cor} If $\alpha$ and $\beta$ are positive numbers such that $\alpha \beta  = {\pi ^2}$, and $w>0$, then
\begin{align}
\frac{\pi}{2\sin(\sqrt{w\alpha})\sinh(\sqrt{w\beta})}=\frac{1}{2w}+\su \left\{\frac{n\alpha(-1)^n}{\sinh(n\alpha)(w+n^2\alpha)}+\frac{n\beta(-1)^n}{\sinh(n\beta)(w-n^2\beta)} \right\}.
\end{align}
\end{cor}
\pf Letting $k=1$ in (\ref{5.7}) gives
\[\su \frac{(-1)^n}{n}\left\{\frac{1}{\sinh(n\beta)}-\frac{1}{\sinh(n\alpha)}\right\}=\frac{\beta-\alpha}{12}.\]
Then, setting $A=B=A_2$ in (\ref{5.17}) yields the desired evaluation.\hfill$\square$\\[5mm]
{\bf Acknowledgments.}  The author expresses his deep gratitude to Professor Jianqiang Zhao for valuable discussions and comments.


{\small
\begin{thebibliography}{99}

\bibitem{A2000}
G.E. Andrews, R. Askey and R. Roy, {\sl Special functions},
Cambridge University Press, 2000: 481-532.

\bibitem{AB2009}
G.E. Andrews and B.C. Berndt, {\sl Ramanujan's lost notebook, Part II}, Springer, New York, 2009, 223-233.

\bibitem{AB2013}
G.E. Andrews and B.C. Berndt, {\sl Ramanujan's lost notebook, Part IV}, Springer, New York, 2013, 265-284.

\bibitem{Au2020}
K.\ C.\ Au, Evaluation of one-dimensional polylogarithmic
integral, with applications to infinite series, arxiv: 2007.03957v1.
Mathematica package available at researchgate.net/publication/342344452.

\bibitem{BBG1994}
D.H. Bailey, J.M. Borwein and R. Girgensohn, {\sl Experimental evaluation of Euler sums},
Exp. Math., 1994, {\bf 3}(1): 17-30.

\bibitem{B1989} B.C. Berndt, {\sl Ramanujan's notebooks, Part II}, Springer-Verlag, New York, 1989.

\bibitem{B2004}
B.C. Berndt, {\sl An unpublished manuscript of Ramanujan on infinite series identities}, J. Ramanujan
Math. Soc., 2004, {\bf 19}: 57--74.

\bibitem{BBV2010}
J. Bl${\rm \ddot{u}}$mlein, D.J. Broadhurst and J.A.M. Vermaseren, {\sl The multiple zeta value data mine}. Comput. Phys. Commun., 2012, {\bf 181}(3): 582-625.

\bibitem{BBG1995}
D. Borwein, J.M. Borwein and R. Girgensohn, {\sl Explicit evaluation of
Euler sums}, Proc. Edinburgh Math., 1995, {\bf 38}: 277-294.

\bibitem{BG1996}
J.M. Borwein and R. Girgensohn, {\sl Evaluation of triple Euler sums}, Electron. J. Combin., 1996: 2-7.

\bibitem{Br2002}
D.M. Bradley, {\sl Series acceleration formulas for Dirichlet series with periodic coefficients}, Ramanujan J., 2002, {\bf 6}: 331-346.

\bibitem{FS1998}
P. Flajolet and B. Salvy, {\sl Euler sums and contour integral representations}, Exp. Math., 1998, {\bf 7}(1): 15--35.

\bibitem{H1992}
M.E. Hoffman, {\sl Multiple harmonic series}, Pacific J. Math., 1992, {\bf 152}: 275-290.

\bibitem{K1996}
K. K$\ddot{\rm o}$lbig, {\sl The polygamma function $\psi(x)$ for $x = 1/4 $ and $x = 3/4$}, J. Comput. Appl. Math., 1996, {\bf 75}: 43-46.

\bibitem{M2014}
I. Mez$\ddot{\rm o}$, {\sl Nonlinear Euler sums}, Pacific J. Math., 2014, {\bf 272}: 201-226.

\bibitem{M1925}
S.L. Malurkar, {\sl On the application of Herr Mellin's integrals to some series}, J. Indian Math. Soc., 1925-1926: 130-138.

\bibitem{N1906}
N. Nielsen, {\sl Handbuch der Theorie der
Gammafunktion and Theorie des Integrallogarithmua
und ueruumdier Transzendenten}, 1906. Reprinted
together as Die Gammafunktion, Chelsea, New York,
1965.

\bibitem{R2012}
S. Ramanujan, {\sl Notebooks (2 volumes)}, Tata Institute of Fundamental Research, Bombay, 1957; seconded, 2012.

\bibitem{R1988}
S. Ramanujan, {\sl The lost notebook and other unpublished papers}, Narosa, New Delhi, 1988.

\bibitem{S1}
R. Sitaramachandrarao, {\sl Ramanujan's formula for $\zeta(2n + 1)$}, Madurai Kamaraj University Technical
Report 4, pp. 70--117.

\bibitem{W2017}
W. Wang and Y. Lyu, {\sl Euler sums and Stirling sums}, J. Number Theory, 2018, {\bf 185}: 160-193.

\bibitem{Xu2017}
C. Xu, {\sl Multiple zeta values and Euler sums}, J. Number Theory, 2017, {\bf 177}: 443-478.

\bibitem{XW2018}
C. Xu and W. Wang, {\sl Explicit formulas of Euler sums via multiple zeta values}, J. Symb. Comput., 2020, {\bf 101}: 109-127.

\bibitem{DZ1994}
D. Zagier, {\sl Values of zeta functions and their applications}. First European Congress
of Mathematics, Volume II, Birkhauser, Boston., 1994, (120): 497-512.

\bibitem{Z2016}
J. Zhao, {\sl Multiple zeta functions, multiple polylogarithms and their special values, Series on Number
Theory and its Applications}, 12, World Scientific Publishing Co. Pte. Ltd., Hackensack, NJ, 2016.

\bibitem{Z2019}
M. Zhao, {\sl On specific log integrals, polylog integrals and alternating Euler sums}, arXiv:1911.12155v5.

\end{thebibliography}}
 \end{document}